\documentclass[11pt,a4paper]{article}
\usepackage{a4wide}
\usepackage[T1]{fontenc}
\usepackage{latexsym}
\usepackage{amsmath}
\usepackage{amsthm}
\usepackage{amssymb,enumerate}
\usepackage[usenames]{color}
\usepackage{graphicx}
\usepackage{fancybox}
\usepackage[normalem]{ulem}
\usepackage{hyperref}

\theoremstyle{plain}
\newtheorem{theorem}{Theorem}[section]
\newtheorem{proposition}[theorem]{Proposition}

\newtheorem{lemma}[theorem]{Lemma}
\newtheorem{prop}[theorem]{Proposition}

\newtheorem*{coro*}{Corollary}
\theoremstyle{definition}
\newtheorem{remark}[theorem]{Remark}
\newtheorem{definition}[theorem]{Definition}

\newtheorem{examples}[theorem]{Examples}
\numberwithin{equation}{section}

\newcommand{\CC}{\mathbb{C}}
\newcommand{\N}{\mathbb{N}}
\newcommand{\R}{\mathbb{R}}
\newcommand{\C}{\mathbb{C}}

\newcommand{\al}{\alpha}
\newcommand{\be}{\beta}

\newcommand{\bL}{{\boldsymbol{L}}}
\newcommand{\m}{{\boldsymbol{m}}}
\newcommand{\bl}{{\boldsymbol{\ell}}}
\newcommand{\M}{{\boldsymbol{M}}}
\newcommand{\NN}{{\boldsymbol{N}}}

\newcommand{\bGa}{{\boldsymbol{\Gamma}}}
\newcommand{\LL}{{\boldsymbol{L}}}
\newcommand{\hM}{\widehat{{\boldsymbol{M}}}}

\newcommand{\bbepsilon}{{\boldsymbol{\varepsilon}}}

\def\a{\alpha}

\def\ga{\gamma}
\def\L{\mathbb{L}}

\newcommand{\lc}{\operatorname{(lc)}}
\newcommand{\dc}{\operatorname{(dc)}}
\newcommand{\mg}{\operatorname{(mg)}}
\newcommand{\snq}{\operatorname{(snq)}}
\newcommand{\sm}{\operatorname{(sm)}}

\definecolor{azulosc}{rgb}{0.2,0.1,0.7}
\definecolor{granate}{rgb}{0.6,0,0.3} 
\definecolor{verdeosc}{RGB}{46, 139, 87} 
\definecolor{rojo}{RGB}{219,0,0}                

\makeatletter
\DeclareRobustCommand\widecheck[1]{{\mathpalette\@widecheck{#1}}}
\def\@widecheck#1#2{%
    \setbox\z@\hbox{\m@th$#1#2$}%
    \setbox\tw@\hbox{\m@th$#1%
       \widehat{%
          \vrule\@width\z@\@height\ht\z@
          \vrule\@height\z@\@width\wd\z@}$}%
    \dp\tw@-\ht\z@
    \@tempdima\ht\z@ \advance\@tempdima2\ht\tw@ \divide\@tempdima\thr@@
    \setbox\tw@\hbox{%
       \raise\@tempdima\hbox{\scalebox{1}[-1]{\lower\@tempdima\box
\tw@}}}%
    {\ooalign{\box\tw@ \cr \box\z@}}}
\makeatother

\begin{document}

\title{Surjectivity of the asymptotic Borel map in Carleman ultraholomorphic classes defined by sequences with shifted moments}
\author{Javier Jim\'enez-Garrido \and Ignacio Miguel-Cantero \and Javier Sanz \and Gerhard Schindl}
\date{\today}

\maketitle

\begin{abstract}
We prove several improved versions of the Borel-Ritt theorem about the surjectivity of the asymptotic Borel mapping in classes of functions with $\M$-uniform asymptotic expansion on an unbounded sector of the Riemann surface of the logarithm. While in previous results the weight sequence $\M$ of positive numbers is supposed to be derivation closed, a much weaker condition is shown to be sufficient to obtain the result in the case of Roumieu classes. Regarding Beurling classes, we are able to slightly improve a classical result of J. Schmets and M. Valdivia and reprove a result of A. Debrouwere, both under derivation closedness. Our new condition also allows us to obtain surjectivity results for Beurling classes in suitably small sectors, but the technique is now adapted from a classical procedure already appearing in the work of V. Thilliez, in its turn inspired by that of J. Chaumat and A.-M. Chollet.
\par\medskip

 \end{abstract}

\section{Introduction}
This paper intends to improve known results on the surjectivity of the asymptotic Borel map in so-called ultraholomorphic Carleman classes in an unbounded sector, both of Roumieu or Beurling type. Such classes of holomorphic functions are defined by restricting the growth of their derivatives, or that of the remainders in their uniform asymptotic expansion, in terms of a given weight sequence $\M=(M_p)_{p\in\N_0}$ of positive real numbers, and the Borel map sends a function into the formal power series providing its expansion. The classical Borel-Ritt-Gevrey theorem of B.~Malgrange and J.-P.~Ramis~\cite{Ramis1} solved the problem for Gevrey asymptotics, for which $\M=(p!^{\alpha})_{p\in\N_0}$, $\alpha>0$. More general situations were covered by J.~Schmets and M.~Valdivia~\cite{SchmetsValdivia00}, V.~Thilliez~\cite{Thilliez03}, A. Debrouwere and the authors~\cite{SanzFlatProxOrder,JimenezSanzSchindlInjectSurject, momentsdebrouwere,DebrouwereBorelRittBeurlingClasses, JimenezSanzSchindlSurjectDC,JimenezMiguelSanzSchindlOptFlat}.
The most satisfactory results are valid for regular weight sequences in the sense of E.~M.~Dyn'kin~\cite{Dynkin80} (see Subsection~\ref{subsectstrregseq} for the definitions). They determine the length of the interval of surjectivity, consisting of those values $\ga>0$ such that the Borel map is surjective when considered on classes defined in the unbounded sector $S_{\gamma}$ bisected by the positive real line and with opening $\pi\gamma$ in the Riemann surface of the logarithm. Such length is precisely $\gamma(\M)$, an index introduced by V. Thilliez~\cite{Thilliez03} for strongly regular sequences and later studied~\cite{JimenezSanzSchindlIndex} for any weight sequence. The first proof of this fact for Roumieu classes, given in~\cite{JimenezSanzSchindlSurjectDC} and resting on the result of A. Debrouwere~\cite{momentsdebrouwere} for a halfplane, had an abstract, functional-analytic flavor. Subsequently, and in the same vein as in a previous work of A. Lastra, S. Malek and the third author~\cite{LastraMalekSanzContinuousRightLaplace} for strongly regular sequences, the construction of optimal flat functions in sectors $S_{\ga}$ with $0<\ga<\ga(\M)$ and for general weight sequences allowed the authors~\cite{JimenezMiguelSanzSchindlOptFlat} to provide a constructive proof by means of formal Borel- and truncated Laplace-like transforms, defined from suitable kernel functions obtained from those optimal flat functions. Moreover, surjectivity comes with local extension operators, linear and continuous right inverses for the Borel map when acting on suitable Banach spaces within our classes.

The main aim of this paper is to put forward a new condition for weight sequences, much weaker than the condition of derivation closedness included in the definition of regular weight sequences, and still allowing for the obtention of Borel-Ritt theorems in our setting in a constructive way. We say $\M=(M_p)_p$ has shifted moments, $\sm$ for short, if there exist $C_0>0$ and $H>1$ such that $$
\log\left(\frac{m_{p+1}}{m_p}\right)\le C_0H^{p+1},\ \ p\in\{0,1,2,\dots\},
$$
where $m_p=M_{p+1}/M_p$. It turns out that, whenever $\ga(\M)>0$, $\sm$ amounts to the equivalence of $\M_{+1}:=(M_{p+1})_p$ and the sequence of Stieltjes moments for a kernel $e(z)=G(1/z)$ defined from an optimal flat function $G$ in the class defined by $\M$. As we will prove, under this weak condition it is possible to adapt the Borel- and truncated Laplace-transforms in order to make our technique work and obtain local extension operators, and so the surjectivity, of the Borel map for Roumieu classes.

Regarding the Beurling case, A. Debrouwere~\cite[Th. 7.4]{momentsdebrouwere} first characterized the surjectivity of the asymptotic Borel map in the right half-plane for regular sequences, and later on he completely solved the problem for non-uniform asymptotics, and provided  global extension operators for $\ga<\ga(\M)$ in the case with uniform estimates, see~\cite{DebrouwereBorelRittBeurlingClasses}. We will present a different technique in order to treat the problem for classes with uniform estimates, following the same ideas as in the Roumieu case~\cite{JimenezSanzSchindlSurjectDC}, which rest on the use of ramified Borel and Laplace integral transforms. In order to do this, we need to prove Theorem~\ref{th.SurjectivityUniformAsymp.plus.dc}, which slightly improves both a result of J. Schmets and M. Valdivia~\cite{SchmetsValdivia00}, Theorem 4.1 in this paper, and the implication $(i)\Rightarrow (iii)$ of the aforementioned result of A. Debrouwere (Theorem 4.4 in this paper).
Finally, the new condition $\sm$ is also valid in order to prove surjectivity for Beurling classes as long as $0<\ga<\ga(\M)$, thanks to a technique of J. Chaumat and A. M. Chollet~\cite{chaucho} already applied by V. Thilliez~\cite[Th. 3.4.1]{Thilliez03} for strongly regular sequences.

The paper is organized as follows. Section~\ref{sectPrelimin} contains the preliminaries, namely: the main information concerning weight sequences, indices or auxiliary functions associated with them, different ultraholomorphic classes and the asymptotic Borel map, and the known result for Roumieu classes by the use of optimal flat functions. In Section~\ref{sect.Improvements.Roumieu} we obtain the new surjectivity results for Roumieu classes defined by weight sequences satisfying $\sm$, while Section~\ref{sect.Improvements.Beurling} contains the improvements for Beurling classes, according to whether $\dc$ or $\sm$ is assumed.

The results included in this paper are part of the Ph.D. dissertation of the second author~\cite[Ch.~3]{PhDDissertationIgnacioMiguel}, defended at the University of Valladolid (Spain) under the advice of the third author.

\section{Preliminaries}\label{sectPrelimin}

\subsection{Weight sequences and their properties}\label{subsectstrregseq}

We set $\N:=\{1,2,...\}$, $\N_{0}:=\N\cup\{0\}$. In what follows, $\M=(M_p)_{p\in\N_0}$ will always stand for a sequence of positive real numbers with $M_0=1$. We define its {\it sequence of quotients}  $\m=(m_p)_{p\in\N_0}$ by
$m_p:=M_{p+1}/M_p$, $p\in \N_0$. The knowledge of $\M$ amounts to that of $\m$, since $M_p=m_0\cdots m_{p-1}$, $p\in\N$. We will denote by small letters the quotients of a sequence given by the corresponding capital letters.
Furthermore, we set $\widehat{\M}:=(p!M_p)_{p \in \N_0}$ and $\widecheck{\M}:=(M_p/p!)_{p \in \N_0}$, and $\widehat{\m}=(\widehat{m}_p)_p$ and $\widecheck{\m}=(\widecheck{m}_p)_p$ will denote the corresponding sequences of quotients: $\widehat{m}_p=(p+1)m_p$, $\widecheck{m}_p=m_p/(p+1)$ for $p\in\N_0$.

We shall use the following conditions on  sequences $\M$:

\begin{itemize}
\item[$\lc$] $\M$ is \emph{log-convex} if $M^2_p \leq M_{p-1}M_{p+1}$, $p \in \N$. Equivalently, $(m_p)_{p\in\N_0}$ is nondecreasing.
\item[$\sm$] $\M$ has \emph{shifted moments} if $\log(m_{p+1}/m_p) \leq C_0H^{p+1}$, $p \in \N_0$, for some $C_0>0$ and $H>1$.
\item[$\dc$] $\M$ is \emph{derivation closed} if $M_{p+1} \leq C_0H^{p+1}M_p$, $p \in \N_0$, for some $C_0>0$ and $H>1$.
\item[$\mg$] $\M$ has \emph{moderate growth} if $M_{p+q} \leq C_0H^{p+q}M_pM_q$, $p,q \in \N_0$, for some $C_0>0$ and $H>1$.
\item[$\snq$] $\M$ is \emph{strongly non-quasianalytic} if $\displaystyle \sum_{q=p}^\infty \frac{1}{(q+1)m_q} \leq  \frac{C}{m_p}$, $p \in \N_0$, for some $C > 0$.
\end{itemize}

\begin{remark}
The form of the inequalities in $\sm$, $\dc$ and $\mg$ admits slight modifications. For example, the constant $C_0$ could be omitted in all of them, or one could use $C_0H^p$ instead of $C_0H^{p+1}$ in $\sm$ or $\dc$.
The choice of one of these alternative expressions can provide some flexibility in our arguments, and we will use one or another without further comment.
\end{remark}

In the classical work of H.~Komatsu~\cite{Komatsu}, the properties $\lc$, $\dc$ and $\mg$ are denoted by $(M.1)$, $(M.2)'$ and $(M.2)$, respectively, while $\snq$ for $\M$ is the same as property
$(M.3)$ for $\widehat{\M}$.

If $\M$ is $\lc$, not only $(m_p)_{p\in\N_0}$ but also $(M_p^{1/p})_{p\in\N}$ is nondecreasing, and $(M_p)^{1/p}\leq m_{p-1}$ for every $p\in\N$; moreover, $\lim_{p\to\infty} (M_p)^{1/p}= \infty$ if and only if $\lim_{p\to\infty} m_p= \infty$. The sequence $\M$ is said to be a \emph{weight sequence} if it satisfies $\lc$ and $m_p \nearrow \infty$ as $p \to \infty$.

According to E.~M.~Dyn'kin~\cite{Dynkin80}, if $\M$ is a weight sequence and satisfies $\dc$, we say $\hM$ is \emph{regular}. Following V.~Thilliez~\cite{Thilliez03}, if $\M$ satisfies $\lc$, $\mg$ and $\snq$, we say $\M$ is \emph{strongly regular}; in this case $\M$ is a weight sequence, and the corresponding $\hM$ is regular (since $\mg$ clearly implies $\dc$).

\begin{remark}\label{rem_propertiesMandstability}
All the previously listed properties are preserved when passing from $\M$ to $\widehat{\M}$. However, only $\sm$, $\dc$ and $\mg$ are generally kept when going from $\M$ to $\widecheck{\M}$.
While these statements are well-known for $\dc$ and $\mg$, the ones for the, up to our knowledge, new condition $\sm$ stem from the inequalities
$$
\log\left(\frac{m_{p+1}}{m_p}\right)\le \log\left(\frac{m_{p+1}}{m_p}\right)+\log\left(\frac{p+2}{p+1}\right)=\log\left(\frac{\widehat{m}_{p+1}}{\widehat{m}_p}\right) \le\log\left(\frac{m_{p+1}}{m_p}\right)+\log(2)
$$
(when applied to $\M$ or to $\widecheck{\M}$). Similarly, $\sm$ (like $\dc$ and $\mg$) is kept when going from $\M$ to $(p!^{\a}M_p)_p$ for any real number $\a$.

Regarding $\sm$ and $\dc$ we have the following.
\end{remark}

\begin{lemma}\label{lemma.dc-implies-sm}
Let $\M$ be a sequence such that $a_0:=\inf_{p\in\N_0}m_p>0$ (in particular, this holds if $\M$ is $\lc$). Then, $\dc$ implies $\sm$.
\end{lemma}
\begin{proof}
Since $m_p\le C_0H^{p+1}$ and $p+1\le 2^{p}$ for every $p\in\N_0$ and some $C_0>0$ and $H>1$, one has
$$
\log\left(\frac{m_{p+1}}{m_p}\right)\le \log(m_{p+1})-\log(a_0)\le\log(C_0/a_0)+(p+2)\log(H)\le C_1H_1^{p+1}
$$
for every $p\in\N_0$, with the choices $C_1=\log(HC_0^{1/2}/a_0^{1/2})>0$ (note that $a_0\le C_0H<C_0H^2$) and $H_1=2$.
\end{proof}

The converse statement in the previous lemma does not hold, as shown by the forthcoming Example~\ref{examp.Sequences}.(iii).

\begin{remark}
Under the same hypothesis of the previous Lemma, one easily checks that the condition
$$
\exists C_1>0,\ H>1\colon \log(m_p)\le C_1H^p,\quad p\in\N_0\eqno{(*)}
$$
implies $\sm$. Indeed, these two conditions are then equivalent, since $\sm$ implies that
$$
\log(m_p)=\log(m_0)+\sum_{k=0}^{p-1}\log\left(\frac{m_{k+1}}{m_k}\right)
\le\log(m_0)+C_0\sum_{k=0}^{p-1}H^{k+1}=\log(m_0)+C_0\frac{H^{p+1}-H}{H-1},
$$
and $(*)$ is readily obtained. Although the condition $(*)$ is simpler, and easier to check in concrete examples, than $\sm$, we have preferred to keep the latter as our standard assumption since it enters our arguments directly, see for example the proof of Proposition~\ref{prop.MomentsEquivM+1}.
\end{remark}

\begin{examples}\label{examp.Sequences}
\begin{itemize}
\item[(i)] The Gevrey sequences $(p!^{\alpha})_{p}$ ($\alpha>0$) are strongly regular, and their perturbations $(p!^{\alpha}\prod_{j=0}^p\log^{\beta}(e+j))_{p}$ ($\alpha>0$, $\beta\in\R$) also are after suitably modifying a finite number of their terms, if necessary. They appear all through the study of formal power series solutions to differential and difference equations.
\item[(ii)] The sequences $\M=(q^{p^\alpha})_{p}$ ($q>1$, $0<\alpha\le 2$) and $\M=(p^{\tau p^\sigma})_{p}$ ($\tau>0$, $1<\sigma<2$) are such that $\hM$ is regular, but they are not strongly regular since they do not satisfy $\mg$. In particular, the $q$-Gevrey sequences $\M=(q^{p^2})_{p}$ ($q>1$) appear in the study of $q$-difference equations.
\item[(iii)] The weight sequences $\M=(q^{p^\alpha})_{p}$ ($q>1$, $\alpha>2$) and $\M=(p^{\tau p^\sigma})_{p}$ ($\tau>0$, $\sigma\ge 2$) do not satisfy $\dc$, so that $\hM$ is not regular, but they still satisfy $\sm$. The sequences of the family $\{(p^{\tau p^\sigma})_{p}\}_{\sigma>1}$ have appeared as the defining sequences for some generalized ultradifferentiable classes ``beyond Gevrey regularity'', deeply studied in a series of papers by S. Pilipovi{\'c}, N. Teofanov and F. Tomi{\'c}~\cite{ptt15,ptt16,ptt20,ptt21,TeofanovTomic} and by some of the authors of this paper and their collaborators~\cite{JimenezLastraSanzRapidGrowth,JimenezNenningSchindl}.
\item[(iv)] The rapidly growing weight sequences $\M=(q^{p^p})_{p}$ ($q>1$; $M_0:=1$) are weight sequences which do not satisfy $\sm$. As it will be seen, they are the only ones in the list to which the results in this paper cannot be applied.
\end{itemize}
\end{examples}

Following H. Komatsu~\cite{Komatsu}, the relation $\M \subset \NN$ between two sequences means that
there are $C,h > 0$ such that $M_p \leq Ch^pN_p$ for all $p \in \N_0$. We say $\M$ and $\NN$ are \emph{equivalent}, denoted  $\M \approx \NN$, if $\M \subset \NN$ and $\NN \subset \M$; equivalently, there exist $h_1,h_2>0$ such that $h_1^{p+1}N_p\le M_p \leq h_2^{p+1}N_p$ for all $p \in \N_0$. It is straightforward that $\dc$ and $\mg$ are stable under equivalence for general sequences, and the same can be deduced for $\snq$ for weight sequences by an indirect method, as this condition characterizes the surjectivity of the Borel mapping in Carleman ultradifferentiable classes, by a result of H.-J. Petzsche~\cite[Cor. 3.2]{Petzsche88} (see~\cite[Cor. 3.14]{JimenezSanzSchindlIndex} for a direct proof of a more general statement about the stability of $\snq$). We prove the result for $\sm$.

\begin{lemma}\label{lemma.sm-stable-equivalence}
The property $\sm$ is preserved under equivalence of sequences.
\end{lemma}

\begin{proof}
Suppose $\M$ satisfies $\sm$, and consider $\LL=(L_p)_p$ such that $\M\approx\LL$. There exists $C>1$ such that $C^{-p-1}L_p\le M_p\le C^{p+1}L_p$, $p\in\N_0$.
Consequently, for all $p\in\N_0$ one has
$$
l_p=\frac{L_{p+1}}{L_p}\le\frac{C^{p+2}C^{p+1}M_{p+1}}{M_p}=C^{2p+3}m_p,\ \ l_p\ge\frac{M_{p+1}}{C^{p+2}C^{p+1}M_p}=C^{-2p-3}m_p.
$$
So, taking into account that $p+1\le 2^p$ for every $p$, we have
\begin{align*}
\log\left(\frac{l_{p+1}}{l_p}\right)&\le \log\left(C^{4p+8}\frac{m_{p+1}}{m_p}\right)
\le(4p+8)\log(C)+C_0H^{p+1}\\
&\le 4\log(C)+2\log(C)2^{p+1}+C_0H^{p+1}\le C_1H_1^{p+1}
\end{align*}
for the choices $C_1=4\log(C)+C_0>0$ and $H_1=\max\{2,H\}>1$, so that $\LL$ also satisfies $\sm$.
\end{proof}

Given a sequence $\M=(M_p)_p$ and $r>0$, we write $\M^r:=(M_p^r)_p$. It is clear that $\M$ satisfies $\sm$ if, and only if, $\M^r$ does for some/every $r>0$.

\subsection{Associated function with a sequence}\label{subsectAssocFunctSeq}

It is very helpful to consider the \emph{associated function} to a sequence $\M$ given by
\begin{equation*}
	h_\M(t):=\inf_{p\in\N_0}M_p t^p, \quad t>0;\ h_\M(0):=0.
\end{equation*}
It is well known that $h_{\M}(t)>0$ for every $t>0$ if, and only if, $\lim_{p\to\infty}M_p^{1/p}=\infty$, which will be a standard assumption when $h_{\M}$ enters our considerations. In particular, this condition is satisfied for weight sequences. From~\cite[Prop.~3.2]{Komatsu} we find that, for a weight sequence $\M$,
\begin{equation}\label{eq.MpfromomegaM}
	M_p=
\sup_{t>0}t^p h_{\M}(1/t),\quad p\in\N_0,
\end{equation}
and it is plain to check that
\begin{equation}\label{eq.expression_h_M}
h_\M(t)=M_{p+1}t^{p+1} \quad\textrm{for }\frac{1}{m_{p+1}}\le t <\frac{1}{m_p},\ p\in\N_0;\quad h_\M(t)=1\quad\textrm{for }t\ge\frac{1}{m_0}.
\end{equation}
Moreover, $h_\M$ is nondecreasing and continuous.

\subsection{Index $\ga(\M)$ for weight sequences}\label{subsectIndexGammaM}

The index $\gamma(\M)$, introduced by V.~Thilliez~\cite[Sect.\ 1.3]{Thilliez03} for strongly regular sequences $\M$, can be defined for $\lc$ sequences, and it may be equivalently expressed by different conditions:
\begin{enumerate}[(i)]
\item A sequence $(c_p)_{p\in\N_0}$ is \emph{almost increasing} if there exists $a>0$ such that for every $p\in\N_0$ we have that $c_p\leq a c_q $ for every $ q\geq p$.
It was proved in~\cite{JimenezSanzSRSPO,JimenezSanzSchindlIndex} that for any weight sequence $\M$ one has
\begin{equation}\label{equa.indice.gammaM.casicrec}
\gamma(\M)=\sup\{\gamma>0:(m_{p}/(p+1)^\gamma)_{p\in\N_0}\hbox{ is almost increasing} \}\in[0,\infty].
\end{equation}
\item Given $\be>0$, we say that $\m$ satisfies the condition $(\gamma_{\be})$  if there exists $A>0$ such that
\begin{equation*}
\sum^\infty_{\ell=p} \frac{1}{(m_\ell)^{1/\be}}\leq \frac{A (p+1) }{(m_p)^{1/\be}},  \qquad p\in\N_0.\tag{$\gamma_\be$}
\end{equation*}
In~\cite{PhDJimenez,JimenezSanzSchindlIndex} it is proved that for a weight sequence $\M$,
\begin{equation}\label{equa.indice.gammaM.gamma_r}
\gamma(\M)=\sup\{\be>0\colon \,\, \m \,\, \text{satisfies } (\gamma_{\be})\,\};\quad \gamma(\M)>\be\iff\m\text{ satisfies }(\gamma_\be).
\end{equation}
\end{enumerate}

If we observe that the condition $\snq$ for $\M$ is precisely $(\gamma_1)$ for $\widehat{\m}$, the sequence of quotients for $\hM$, and that $\ga(\hM)=\ga(\M)+1$ (this is clear from~\eqref{equa.indice.gammaM.casicrec}), we deduce from the second statement in~\eqref{equa.indice.gammaM.gamma_r} that
$\M$ satisfies $\snq$ if, and only if, $\ga(\M)>0$.

Let $\M$ and $\LL$ be sequences, and $\m=(m_p)_p$ and $\bl=(\ell_p)_p$ their respective sequences of quotients. We write $\m\simeq\bl$ whenever there exists $c>0$ such that $c^{-1}\ell_p\le m_p\le c\ell_p$, $p\in\N_0$. It is clear that $\m\simeq\bl$ implies $\M\approx\L$.

We recall the following result for later use.

\begin{lemma}[\cite{JimenezSanzSchindlIndex}, Remark 3.15]\label{lemma.gammaMgreaterthan1}
	For an arbitrary sequence $\M$ such that $\gamma(\M)>1$,
	there exists a weight sequence $\LL$ such that $\m\simeq\widehat{\ell}$, and so $\widehat{\LL}\approx\M$ and $\gamma(\widehat{\LL})=\gamma(\M)$.
\end{lemma}

For $\al>0$ we set
$\bGa_{\al}:=(\Gamma(1+\al p))_{p\in\N_0}$, where $\Gamma$ is Euler's Gamma function; its sequence of quotients will be denoted by $\overline{\boldsymbol{g}}^{\al}$.
It is well known that $\bGa_{\al}\approx(p!^{\al})_{p\in\N_0}$, and a straightforward verification shows that for any sequence $\M$ and for every $\al>0$ one has
\begin{align} \gamma((p!^{\al}M_p)_{p\in\N_0})&=\gamma((\Gamma(1+\al p)M_p)_{p\in\N_0})= \gamma(\M)+\al,\label{equa.gamma.producto}\\
\gamma((M_p/p!^{\al})_{p\in\N_0})&=\gamma(M_p/(\Gamma(1+\al p))_{p\in\N_0})= \gamma(\M)-\al.\label{equa.gamma.cociente}
\end{align}

\subsection{Asymptotic expansions, ultraholomorphic classes and the asymptotic Borel map}\label{subsectCarlemanclasses}

In what follows $\mathcal{R}$ stands for the Riemann surface of the logarithm, and $\C[[z]]$ is the space of formal power series in $z$ with complex coefficients.

For $\gamma>0$, we consider unbounded sectors bisected by direction 0,
$$S_{\gamma}:=\{z\in\mathcal{R}:|\hbox{arg}(z)|<\frac{\gamma\,\pi}{2}\}$$
or, in general, unbounded sectors with bisecting direction $d\in\R$ and opening $\ga\,\pi$,
$$S(d,\ga):=\{z\in\mathcal{R}:|\hbox{arg}(z)-d|<\frac{\ga\,\pi}{2}\}.$$

A sector $T$ is said to be a \emph{proper subsector} of a sector $S$ if $\overline{T}\subset S$ (where the closure of $T$ is taken in $\mathcal{R}$, and so the vertex of the sector is not under consideration).

In this paragraph $S$ is an unbounded sector and $\M$ a sequence. We start by recalling the concept of uniform asymptotic expansion.

We say a holomorphic function $f\colon S\to\C$ admits $\widehat{f}=\sum_{n\ge 0}a_nz^n\in\C[[z]]$ as its \emph{uniform $\M$-asymptotic expansion in $S$ (of type $1/h$ for some $h>0$)} if there exists $C>0$ such that for every $p\in\N_0$, one has
\begin{equation}\left|f(z)-\sum_{n=0}^{p-1}a_nz^n \right|\le Ch^pM_{p}|z|^p,\qquad z\in S.\label{desarasintunifo}
\end{equation}
In this case we write $f\sim_{\M,h}^u\widehat{f}$ in $S$, and $\widetilde{\mathcal{A}}^u_{\M,h}(S)$ denotes the space of functions admitting uniform $\M$-asymptotic expansion of type $1/h$ in $S$, endowed with the norm
\begin{equation}\label{eq.NormUniformAsymptFixedType}
\left\|f\right\|_{\M,h,\overset{\sim}{u}}:=\sup_{z\in S,p\in\N_{0}}\frac{|f(z)-\sum_{k=0}^{p-1}a_kz^k|}{h^{p}M_{p}|z|^p},
\end{equation}
which makes it a Banach space. $\widetilde{\mathcal{A}}^u_{\{\M\}}(S)$ stands for the $(LB)$ space of functions admitting a uniform $\{\M\}$-asymptotic expansion in $S$, obtained as the union of the previous classes when $h$ runs over $(0,\infty)$. When the type needs not be specified, we simply write $f\sim_{\{\M\}}^u\widehat{f}$ in $S$.
It is also interesting to consider the space of Beurling-type $\widetilde{\mathcal{A}}^u_{(\M)}(S):=\cap_{h>0}\widetilde{\mathcal{A}}^u_{\M,h}(S)$, which becomes a Fr\'echet space when endowed with the topology generated by the family of (semi)norms $(\left\|\cdot\right\|_{\M,h,\overset{\sim}{u}})_{h>0}$.

Note that, taking $p=0$ in~\eqref{desarasintunifo}, we deduce that every function in $\widetilde{\mathcal{A}}^u_{\{\M\}}(S)$ or $\widetilde{\mathcal{A}}^u_{(\M)}(S)$ is a bounded function.

Finally, we define for every $h>0$ the class $\mathcal{A}_{\M,h}(S)$ consisting of the holomorphic functions $f$ in $S$ such that
$$
\left\|f\right\|_{\M,h}:=\sup_{z\in S,p\in\N_{0}}\frac{|f^{(p)}(z)|}{h^{p}M_{p}}<\infty.
$$
($\mathcal{A}_{\M,h}(S),\left\|\,\cdot\, \right\|_{\M,h}$) is a Banach space, and $\mathcal{A}_{\{\M\}}(S):=\cup_{h>0}\mathcal{A}_{\M,h}(S)$ is called a \emph{Carleman-Roumieu ultraholomorphic class} in the sector $S$, whose natural inductive topology makes it an $(LB)$ space. On the other hand, $\mathcal{A}_{(\M)}(S):=\cap_{h>0}\mathcal{A}_{\M,h}(S)$ is a \emph{Carleman-Beurling ultraholomorphic class} in the sector $S$, whose Fr\'echet space structure is given by the family of seminorms $(\left\|\cdot\right\|_{\M,h})_{h>0}$.\par

We warn the reader that these notations, while the same as in the papers~\cite{JimenezSanzSchindlSurjectDC,JimenezMiguelSanzSchindlOptFlat}, do not agree with the ones used in~\cite{SanzFlatProxOrder,JimenezSanzSchindlInjectSurject}, where
$\widetilde{\mathcal{A}}^u_{\{\M\}}(S)$ was denoted by
$\widetilde{\mathcal{A}}^u_{\M}(S)$,
$\mathcal{A}_{\M,h}(S)$ by $\mathcal{A}_{\widecheck{\M},h}(S)$, and $\mathcal{A}_{\{\M\}}(S)$ by $\mathcal{A}_{\widecheck{\M}}(S)$, respectively.

When a statement is valid for both Roumieu and Beurling classes, we will use the notation $\mathcal{A}_{[\M]}(S)$, $\widetilde{\mathcal{A}}^u_{[\M]}(S)$ and so on (substituting every square bracket by either of them, curly brackets or parentheses, but the same all through the statement).

If $\M$ is $\lc$, the spaces $\mathcal{A}_{[\M]}(S)$ and $\widetilde{\mathcal{A}}^u_{[\M]}(S)$
are algebras, and if $\M$ is $\dc$ they are stable under taking derivatives.
Moreover, if $\M\approx\bL$ the corresponding classes coincide.

Since the derivatives of $f\in\mathcal{A}_{\M,h}(S)$ are Lipschitz, for every $p\in\N_{0}$ one may define
\begin{equation}\label{eq.deriv.at.0.def}
f^{(p)}(0):=\lim_{z\in S,z\to0 }f^{(p)}(z)\in\C.
\end{equation}

As a consequence of Taylor's formula and Cauchy's integral formula for the derivatives, there is a close relation between Carleman ultraholomorphic classes and the concept of asymptotic expansion (for a similar proof see~\cite[Prop. 8]{balserutx}).

\begin{prop}\label{propcotaderidesaasin-RoumieuBeurling}
Let $\M$ be a sequence and $S$ be a sector.
\begin{enumerate}[(i)]
\item If $f\in\mathcal{A}_{\hM,h}(S)$ then $f$ admits $\widehat{f}:=\sum_{p\in\N_0}\frac{1}{p!}f^{(p)}(0)z^p$ as its uniform $\M$-asymptotic expansion in $S$ of type $1/h$, where $(f^{(p)}(0))_{p\in\N_0}$ is given by \eqref{eq.deriv.at.0.def}. Moreover, $\|f\|_{\M,h,\overset{\sim}{u}}\le \|f\|_{\hM,h}$, and so the identity map $\mathcal{A}_{\hM,h}(S)\hookrightarrow \widetilde{\mathcal{A}}^u_{\M,h}(S)$ is continuous. Consequently, we also have that
$\mathcal{A}_{[\hM]}(S)\subseteq \widetilde{\mathcal{A}}^u_{[\M]}(S)$
and $\mathcal{A}_{[\hM]}(S)\hookrightarrow \widetilde{\mathcal{A}}^u_{[\M]}(S)$ is continuous.

\item If $S$ is unbounded and $T$ is a proper subsector of $S$, then there exists a constant $c=c(T,S)>0$ such that the restriction to $T$, $f|_T$, of functions $f$ defined on $S$ and admitting a uniform $\M$-asymptotic expansion in $S$ of type $1/h>0$, belongs to $\mathcal{A}_{\hM,ch}(T)$, and $\|f|_T\|_{\hM,ch}\le \|f\|_{\M,h,\overset{\sim}{u}}$. So, the restriction map from $\widetilde{\mathcal{A}}^u_{\M,h}(S)$ to $\mathcal{A}_{\hM,ch}(T)$ is continuous, and it is also continuous from $\widetilde{\mathcal{A}}^u_{[\M]}(S)$ to $\mathcal{A}_{[\hM]}(T)$.
\end{enumerate}
\end{prop}

One may similarly define classes of formal power series
\begin{equation}\label{eq.defBanachFormalPowerSeries}
\C[[z]]_{\M,h}=\Big\{\widehat{f}=\sum_{p=0}^\infty a_pz^p\in\C[[z]]:\, \left|\,\widehat{f} \,\right|_{\M,h}:=\sup_{p\in\N_{0}}\displaystyle \frac{|a_{p}|}{h^{p}M_{p}}<\infty\Big\}.
\end{equation}%
$(\C[[z]]_{\M,h},\left| \,\cdot\,  \right|_{\M,h})$ is a Banach space and we set $\C[[z]]_{\{\M\}}:=\cup_{h>0}\C[[z]]_{\M,h}$, the Carleman-Roumieu-type weighted space of formal power series, and $\C[[z]]_{(\M)}:=\cap_{h>0}\C[[z]]_{\M,h}$, the Carleman-Beurling-type weighted space.
$\C[[z]]_{\{\M\}}$ is an $(LB)$ space and $\C[[z]]_{(\M)}$ is a Fr\'echet space, when endowed with their natural locally convex topologies.

The \emph{asymptotic Borel map}  $\widetilde{\mathcal{B}}$
sends a function $f\in\widetilde{\mathcal{A}}^u_{\M,h}(S)$
into its $\M$-asymptotic expansion $\widehat{f}\in\C[[z]]_{\M,h}$.
By Proposition~\ref{propcotaderidesaasin-RoumieuBeurling}.(i)
the asymptotic Borel map may be defined from $\widetilde{\mathcal{A}}^u_{[\M]}(S)$ or $\mathcal{A}_{[\hM]}(S)$ into  $\C[[z]]_{[\M]}$ (with the aforementioned meaning), and from $\mathcal{A}_{\hM,h}(S)$ into $\C[[z]]_{\M,h}$, and it is continuous when considered between the corresponding $(LB)$, Fr\'echet or Banach spaces.

If $\M$ is $\lc$, $\widetilde{\mathcal{B}}$ is a homomorphism of algebras; if $\M$ is also $\dc$, differentiation commutes with $\widetilde{\mathcal{B}}$. Finally, $\M\approx\bL$ implies $\C[[z]]_{[\M]}=\C[[z]]_{[\bL]}$, and the corresponding Borel maps are in all cases identical.

We will focus on the surjectivity of the Borel map in unbounded sectors $S_{\gamma}$ bisected by direction 0, as this problem is invariant under rotation. So, we define
\begin{align*}
S_{[\hM]}:=&\{\gamma>0\colon \quad \widetilde{\mathcal{B}}:\mathcal{A}_{[\hM]}(S_\gamma)\longrightarrow \C[[z]]_{[\M]} \text{ is surjective}\} ,\\
\widetilde{S}^u_{[\M]}:=&\{\gamma>0\colon \quad\widetilde{\mathcal{B}}:\widetilde{\mathcal{A}}^u_{[\M]}(S_\gamma)\longrightarrow \C[[z]]_{[\M]} \text{ is surjective}\}.
\end{align*}

The \emph{surjectivity intervals} $S_{[\hM]}$ and $\widetilde{S}^u_{[\M]}$ are either empty or left-open intervals having $0$ as endpoint. By~Proposition~\ref{propcotaderidesaasin-RoumieuBeurling},
we see that
\begin{align}
(\widetilde{S}^u_{[\M]})^{\circ}\subseteq S_{[\hM]}\subseteq \widetilde{S}^u_{[\M]},
\label{equaContentionSurjectIntervals}
\end{align}
where $I^{\circ}$ is the interior of $I$. The determination of these intervals is closely related to the existence of right inverses for the asymptotic Borel map, i. e., linear and continuous operators $T$ such that $\widetilde{\mathcal{B}}\circ T$ is the identity map on a class of formal power series. They are naturally called \emph{extension operators}, and they can be \emph{global}, defined from  $\C[[z]]_{[\M]}$ onto $\widetilde{\mathcal{A}}^u_{[\M]}(S)$ or $\mathcal{A}_{[\hM]}(S)$ (with their respective $(LB)$ or Fr\'echet space structures), or \emph{local}, at the level of Banach spaces, defined from $\C[[z]]_{\M,h}$ into some $\widetilde{\mathcal{A}}^u_{\M,h'}(S)$ or $\mathcal{A}_{\hM,h'}(S)$ for suitable $h'$ depending on $h$. In this latter case, it is common that a scaling of the type occurs, that is, $h'=ch$ for a universal constant $c>0$ independent from $h$.

\subsection{Optimal flat functions and surjectivity of the Borel map in the Roumieu case and under $\dc$}\label{sectFlatFunctions}

The following result for the Roumieu case, already hinted in the work of V. Thilliez~\cite[Subsect. 3.3]{Thilliez03} and resting on a result of H.-J. Petzsche~\cite[Th. 3.5]{Petzsche88}, appeared, in a slightly different form, in~\cite[Lemma~4.5]{JimenezSanzSchindlInjectSurject}. Since the result of Petzsche is equally valid for the Beurling case~\cite[Th. 3.4]{Petzsche88}, one can state the following.

\begin{lemma}\label{lemma.SurjectImpliesgammaMpositive}
Let $\M$ be a weight sequence. If $\widetilde{S}^u_{[\M]}\neq\emptyset$, then $\M$ satisfies $\snq$ or, equivalently, $\ga(\M) >0$.
\end{lemma}

Regarding the precise determination of the surjectivity intervals, the first seminal results appeared in a work of J. Schmets and M. Valdivia~\cite{SchmetsValdivia00}, whose results prove that $$(0,\lceil\ga(\M)\rceil -1)\subset S_{[\hM]},$$
where $\lceil x\rceil$ is the least integer greater than or equal to a real number $x$. Moreover, for such openings surjectivity comes with local extension operators with scaling of the type, and with global extension operators in the Beurling case, while global extension operators in the Roumieu case need the extra condition $(\beta_2)$ of H.-J. Petzsche~\cite{Petzsche88}.
In the case of strongly regular sequences,  V.Thilliez~\cite{Thilliez03} showed that $(0,\ga(\M))\subset S_{[\hM]}$, again with local extension operators with scaling of the type. Several improvements followed in the Roumieu case~\cite{SanzFlatProxOrder,JimenezSanzSchindlInjectSurject,momentsdebrouwere,JimenezSanzSchindlSurjectDC}, trying firstly to determine the surjectivity intervals, or at least their length, for (certain classes of) strongly regular sequences, and afterwards trying to weaken the condition of moderate growth.
These efforts have led to the following precise statement~\cite[Th.~3.7]{JimenezSanzSchindlSurjectDC} under the condition $\dc$. It shows that the length of the surjectivity intervals is precisely given by $\ga(\M)$.

\begin{theorem}\label{teorSurject.dc}
Let $\hM$ be a regular sequence such that $\gamma(\M)>0$. Then,
$$(0,\gamma(\M))\subseteq S_{\{\hM\}}\subseteq \widetilde{S}^u_{\{\M\}}\subseteq
(0,\gamma(\M)].
$$
In particular, if $\gamma(\M)= \infty$, we have that $S_{\{\hM\}}= \widetilde{S}^u_{\{\M\}}= (0,\infty)$.
\end{theorem}

Moreover, it turns out that the surjectivity of the Borel map in Carleman-Roumieu classes is closely related to the existence of optimal flat functions.

\begin{definition}\label{optimalflatdef}
Let $\M$ be a weight sequence, $S$ an unbounded sector bisected by direction $d=0$, i.e., by the positive real line $(0,+\infty)\subset\mathcal{R}$. A holomorphic function $G\colon S\to\C$ is called an \emph{optimal $\{\M\}$-flat function} in $S$ if:
	\begin{itemize}
		\item[$(i)$] There exist $K_1,K_2>0$ such that for all $x>0$,
		\begin{equation}\label{optimalflatleft}
			K_1h_{\M}(K_2x)\le G(x).
		\end{equation}
		
		\item[$(ii)$] There exist $K_3,K_4>0$ such that for all $z\in S$, one has
		\begin{equation}\label{optimalflatright}
			|G(z)|\le K_3h_{\M}(K_4|z|).
		\end{equation}
	\end{itemize}
\end{definition}

Observe that $G(x)>0$ for $x>0$, and so $G(\overline{z})=\overline{G(z)}$, $z\in S$.
The estimates in~\eqref{optimalflatright} can be rewritten as
$$|G(z)|\le K_3K_4^pM_p|z|^p,\qquad p\in\N_0,\ z\in S,
$$
what exactly means that $G\in\widetilde{\mathcal{A}}^u_{\{\M\}}(S)$ and it is \emph{$\{\M\}$-flat}, i.e., its asymptotic expansion is given by the null series. The inequality imposed in~\eqref{optimalflatleft} makes the function optimal in a sense.

The following result is crucial for our purposes.

\begin{prop}(\cite{JimenezMiguelSanzSchindlOptFlat}, Prop. 3.10)\label{prop.ExistOptFlatFunct}
Let $\M$ be a weight sequence with $\ga(\M)>0$. Then, for any $0<\gamma<\ga(\M)$ there exists an optimal $\{\M\}$-flat function in $S_\ga$.
\end{prop}

By means of an optimal flat function, one can obtain local extension operators with scaling of the type for Carleman-Roumieu ultraholomorphic classes defined by regular sequences. We briefly describe the procedure.

If $G$ is an optimal $\{\M\}$-flat function in $\widetilde{\mathcal{A}}^u_{\{\M\}}(S)$, we define the kernel function $e\colon S\to\C$ given by
\begin{equation*}
e(z):=G\left(\frac{1}{z}\right),\quad z\in S.
\end{equation*}
According to Definition~\ref{optimalflatdef}, $e(x)>0$ for all $x>0$, and there exist $K_1,K_2,K_3,K_4>0$ such that
\begin{equation}\label{eq.Bounds_e_sector}
K_1h_{\M}\left(\frac{K_2}{x}\right)\le e(x), \quad x>0, \hspace{1cm}\text{and}\hspace{1cm} |e(z)|\le K_3h_{\M}\left(\frac{K_4}{|z|}\right),\quad z\in S.
\end{equation}
For every $p\in\mathbb{N}_0$ we define the \emph{$p$-th moment} of the function $e(z)$, given by
\begin{equation*}
\mu(p):=\int_0^\infty t^{p}e(t)\,dt.
\end{equation*}

By carefully inspecting the proof of~\cite[Prop. 3.11]{JimenezMiguelSanzSchindlOptFlat}, where only one implication was mentioned, one can state the following equivalence.

\begin{prop}\label{prop.MomentsEquivM}
Suppose $\M$ is a weight sequence with $\ga(\M)>0$, and $G$ is an optimal $\{\M\}$-flat function in $\widetilde{\mathcal{A}}^u_{\{\M\}}(S)$ for some unbounded sector $S$ bisected by the positive real line. Consider the sequence of moments $\boldsymbol{\mu}:=(\mu(p))_{p\in\N_0}$ associated with the kernel function $e(z)=G(1/z)$. Then, $\M$ satisfies $\dc$ if, and only if,
$\M$ and $\boldsymbol{\mu}$ are equivalent.
\end{prop}

We can already recall the following main result for regular sequences. The previous proposition provides the key fact for the implication $(ii)\Rightarrow(iii)$.

\begin{theorem}(\cite{JimenezMiguelSanzSchindlOptFlat}, Theorem 3.12)\label{teor.surjectRoumieuwithDC}
Let $\hM$ be a regular sequence (that is, $\M$ is a weight sequence and satisfies $\dc$) with $\ga(\M)>0$, and let $\ga>0$ be given. Each of the following statements implies the next one:
\begin{itemize}
\item[(i)] $\ga<\ga(\M)$.
\item[(ii)] There exists an optimal $\{\M\}$-flat function in $\widetilde{\mathcal{A}}^u_{\{\M\}}(S_{\ga})$.
\item[(iii)] There exists $c>0$ such that for every $h>0$ there exists an extension operator from $\C[[z]]_{\M,h}$ into $\widetilde{\mathcal{A}}^u_{\M,ch}(S_{\ga})$, i. e., a linear continuous map $T_{\M,h}\colon\C[[z]]_{\M,h}\to \widetilde{\mathcal{A}}^u_{\M,ch}(S_{\ga})$ such that $\widetilde{\mathcal{B}}\circ T_{\M,h}$ is the identity map in $\C[[z]]_{\M,h}$.
\item[(iv)] The Borel map $\widetilde{\mathcal{B}}\colon \widetilde{\mathcal{A}}^u_{\{\M\}}(S_{\ga})\to\C[[z]]_{\{\M\}}$ is surjective. In other words, $(0,\ga]\subset\widetilde{S}_{\{\M\}}^u$.
\item[(v)] $(0,\ga)\subset S_{\{\hM\}}$.
\item[(vi)] $\ga\le\ga(\M)$.
\end{itemize}
\end{theorem}

We note that the condition $\dc$ is only used in the implications $(ii)\Rightarrow (iii)$ and $(v)\Rightarrow(vi)$.

\begin{remark}\label{rem.extension-operators-bounded-deriv}
Theorem~\ref{teor.surjectRoumieuwithDC} and Proposition~\ref{propcotaderidesaasin-RoumieuBeurling}.$(ii)$ together guarantee that for every $\ga\in(0,\ga(\M))$ there exists $a>0$ such that for every $h>0$ there exists a local extension operator from $\C[[z]]_{\M,h}$ into $\mathcal{A}_{\hM,ah}(S_{\ga})$.
\end{remark}

\section{Improved Borel-Ritt-Gevrey theorem in the Roumieu case}\label{sect.Improvements.Roumieu}

In this section we present our first main result, where the condition $\dc$ is changed into the much weaker condition $\sm$. First, we characterize this new condition in terms of the equivalence of the sequence of moments for a kernel function generated from an optimal flat function and the sequence $\M$ shifted by one. This motivates the name given to the condition $\sm$.

\begin{proposition}\label{prop.MomentsEquivM+1}
	Suppose $\M$ is a weight sequence with $\ga(\M)>0$, and $G$ is an optimal $\{\M\}$-flat function in $\widetilde{\mathcal{A}}^u_{\{\M\}}(S)$, where $S$ is an unbounded sector bisected by the positive real line. Consider the sequence of moments $\boldsymbol{\mu}:=(\mu(p))_{p\in\N_0}$ associated with the kernel function $e(z)=G(1/z)$. Then, $\M$ satisfies $\sm$ if, and only if, $\M_{+1}:=(M_{p+1})_{p\in\N_0}$ and $\boldsymbol{\mu}$ are equivalent.
\end{proposition}

\begin{proof}
Suppose $\M$ satisfies $\sm$.
On the one hand, because of the right-hand inequalities in~\eqref{eq.Bounds_e_sector} and the definition of $h_{\M}$, we have
\begin{equation*}
t^pe(t)\le K_3K_4^pM_p, \ \ p\in\N_0,\ t>0.
\end{equation*}
So, we may write
	\begin{align*}
		\mu(p)&=\int_0^{K_4m_p} t^p e(t)\,dt+\int_{K_4m_p}^{K_4m_{p+1}} t^p e(t)\,dt
		+\int_{K_4m_{p+1}}^\infty \frac{1}{t^2}t^{p+2}e(t)\,dt\\
		&\le K_3 K_4 m_p K_4^pM_p+ K_3 \int_{K_4m_p}^{K_4m_{p+1}} t^p h_\M\left(\frac{K_4}{t}\right)\,dt + K_3K_4^{p+2}M_{p+2}\frac{1}{K_4m_{p+1}}\\
		&= 2K_3 K_4^{p+1} M_{p+1}
		+K_3K_4^{p+1}M_{p+1}\log\left(\frac{m_{p+1}}{m_p}\right),
	\end{align*}
where in the last equality we have used~\eqref{eq.expression_h_M}.
Since there exists $H>1$ such that $\log(m_{p+1}/m_p)\leq H^{p+1}$ for every $p$, we get
	\begin{equation*}
		\mu(p)\le K_3 K_4^{p+1} M_{p+1}(2+H^{p+1})
		\leq  K_3K_4(2+H)(K_4H)^{p}M_{p+1},
	\end{equation*}
	and so $\boldsymbol{\mu}\subset \M_{+1}$.

On the other hand, by the left-hand inequalities in~\eqref{eq.Bounds_e_sector}, for every $p\in\N_0$ we may estimate
	$$
	\mu(p)\ge \int_0^s t^p e(t)\,dt \ge K_1 \int_0^s t^p h_{\M}\left(\frac{K_2}{t}\right)\,dt\ge
	K_1 h_{\M}\left(\frac{K_2}{s}\right)\frac{s^{p+1}}{p+1}.
	$$
	Then, by~\eqref{eq.MpfromomegaM} we deduce that
	$$
	\mu(p)\ge \frac{K_1}{p+1} \sup_{s>0}h_{\M}\left(\frac{K_2}{s}\right)s^{p+1}=
	\frac{K_1}{p+1} K_2^{p+1}M_{p+1}\ge
	K_1K_2\left(\frac{K_2}{2}\right)^{p}M_{p+1},
	$$
	and so $\M_{+1}\subset\boldsymbol{\mu}$, as desired.
	
Conversely, suppose $\M_{+1}\approx\boldsymbol{\mu}$. In particular, there exist $C,h>0$ such that $\mu(p)\le Ch^{p}M_{p+1}$ for $ p\in\N_0$. By the left-hand inequalities in~\eqref{eq.Bounds_e_sector}, we may estimate
	$$
	\mu(p)\ge \int_{K_2m_p}^{K_2m_{p+1}} t^p e(t)\,dt \ge K_1 \int_{K_2m_p}^{K_2m_{p+1}} t^p h_{\M}\left(\frac{K_2}{t}\right)\,dt=K_1K_2^{p+1}M_{p+1} \log\left(\frac{m_{p+1}}{m_p}\right),\ p\in\N_0.
	$$
Therefore,	
$$
	K_1K_2^{p+1}M_{p+1}\log\left(\frac{m_{p+1}}{m_p}\right)\leq Ch^{p}M_{p+1},\quad p\in\N_0,
	$$
and $\M$ satisfies $\sm$.
\end{proof}

\begin{remark} It should be noted that, despite its appearance, the condition $\sm$ is the appropriate condition when working with the tools of Borel $\M$-summability, since the shift arises naturally. More precisely, given a kernel $e(z)$ as before, for a function $f$ belonging to a suitable class we can define the $e-$Laplace transform and the corresponding $e-$Gamma moment function by
$$\mathcal{L}_e (f) (z)=\int_{0}^\infty e(u/z) f(u) du, \qquad m_{e} (p) = \int_{0}^\infty u^{p-1} e(u) du,\quad  p\geq 1.$$
So we have that $\mathcal{L}_e (u^p) (z)= m_e(p+1) z^{p+1}$ for all $p\in\N_0$. Consequently, setting $m_e(0)=1$, these kind of summability techniques can be applied whenever $\M$ and $(m_{e}(p))_{p\in\N_0}$ are equivalent, that is, if and only if $\M$ satisfies $\sm$, since $\mu(p)= m_{e}(p+1)$ for all $p\in\N_0$. Observe that this same shift occurs for the classical Laplace transform given by
$$\mathcal{L}(f) (t)=\int_{0}^\infty e^{-tu} f(u) du,$$
for which  $\mathcal{L}_e (u^p) (t)= \Gamma(p+1)/t^{p+1} =p!/t^{p+1}$ for all $p\in\N_0$.
\end{remark}

We can already state the following main result, whose proof is an adaptation of the one for Theorem~\ref{teor.surjectRoumieuwithDC}. Regrettably, we are not able to deduce $\ga\le\ga(\M)$ from the surjectivity of the Borel map in classes on sectors $S_{\ga}$ under this weaker condition $\sm$.

\begin{theorem}\label{Theorem-surjectivity-logcondition}
	Let $\M$ be a weight sequence satisfying $\sm$ and
with $\ga(\M)>0$, and let $\ga>0$ be given.
Then, each of the following statements implies the next one:
	\begin{itemize}
		\item[(i)] $\ga<\ga(\M)$.
		\item[(ii)] There exists $c>0$ such that for every $h>0$ there exists an extension operator from $\C[[z]]_{\M,h}$ into $\widetilde{\mathcal{A}}^u_{\M,ch}(S_{\ga})$.
		\item[(iii)] The Borel map $\widetilde{\mathcal{B}}\colon \widetilde{\mathcal{A}}^u_{\{\M\}}(S_{\ga})\to\C[[z]]_{\{\M\}}$ is surjective. In other words, $(0,\ga]\subset\widetilde{S}_{\{\M\}}^u$.
		\item[(iv)] $(0,\ga)\subset S_{\{\hM\}}$.
	\end{itemize}
In particular, one has $(0,\ga(\M))\subset S_{\{\hM\}}\subset \widetilde{S}_{\{\M\}}^u$.
\end{theorem}
\begin{proof}
	$(i)\Rightarrow (ii)$ By Proposition~\ref{prop.ExistOptFlatFunct}, valid for any weight sequence $\M$ with $\ga(\M)>0$, we can consider an optimal $\{\M\}$-flat function $G$ in $\widetilde{\mathcal{A}}^u_{\{\M\}}(S_{\ga})$. Let $(\mu(p))_{p\in\N_0}$ be the sequence of moments associated with the function $e(z)=G(1/z)$. Given $h>0$ and $\widehat{f}=\sum_{p=0}^\infty a_pz^p\in\C[[z]]_{\M,h}$, by the definition of the norm in $\C[[z]]_{\M,h}$ (see~\eqref{eq.defBanachFormalPowerSeries}), we have
	\begin{equation*}
		|a_p|\le |\widehat{f}|_{\M,h}h^{p}M_{p},\quad p\in\N_0.
	\end{equation*}
	Because of Proposition~\ref{prop.MomentsEquivM+1}, there exist $h_1,h_2>0$ such that
\begin{equation}\label{eq.boundsMmomentslog}
h_1^{p+1}M_{p+1}\le\mu(p)\le h_2^{p+1}M_{p+1},\ \ p\in\N_0.
\end{equation}
So, we deduce that
\begin{equation}\label{eq.boundsCoeffBorelTransf}
		\left|\frac{a_{p+1}}{\mu(p)}\right|\le |\widehat{f}|_{\M,h}\left(\frac{h}{h_1}\right)^{p+1},\quad p\in\N_0.
	\end{equation}
	Hence, the formal Borel-like transform of $\widehat{f}-a_0$, defined as
	\begin{equation*}
		\widehat{g}=\sum_{p=0}^\infty\frac{a_{p+1}}{\mu(p)}z^p,
	\end{equation*}
	is convergent in the disc $D(0,R)$ for $R=h_1/h>0$, and it defines a holomorphic function $g$ there. Choose $R_0:=h_1/(2h)<R$, and define
	\begin{equation*}
I_{\M,h}(\widehat{f}\,)(z):=\int_{0}^{R_0}e\left(\frac{u}{z}\right)g(u)\,du,\qquad z\in S_{\ga},
	\end{equation*}
	which is a truncated Laplace-like transform of $g$ with kernel $e$. By Leibniz's theorem for parametric integrals and the properties of $e$, this function is holomorphic in $S_{\ga}$. We will prove that
	$I_{\M,h}(\widehat{f}\,)\sim^u_{\{\M\}}\hat{f}-a_0$ uniformly in $S_{\ga}$.
	
Let $p\in\N$ and $z\in S_{\ga}$. We have
	\begin{align*}
		I_{\M,h}(\widehat{f}\,)(z)-\sum_{n=1}^{p-1}a_nz^n &= I_{\M,h}(\widehat{f}\,)(z)-\sum_{n=1}^{p-1}\frac{a_n}{\mu(n-1)}\mu(n-1)z^n\\
		&= \int_{0}^{R_0}e \left(\frac{u}{z}\right) \sum_{n=1}^{\infty}\frac{a_{n}}{\mu(n-1)}u^{n-1}\,du -\sum_{n=1}^{p-1}\frac{a_n}{\mu(n-1)}\int_{0}^{\infty}v^{n-1}e(v)\,dv\, z^n.
	\end{align*}
A change of variable $u=zv$ in the last integral, Cauchy's residue theorem and the right-hand estimates in~(\ref{eq.Bounds_e_sector}) allow us to rotate the path of integration and obtain
	$$
	z^n\int_{0}^{\infty}v^{n-1}e(v)dv= \int_{0}^{\infty}u^{n-1}e\left(\frac{u}{z}\right)\,du.
	$$
	So, the preceding difference can be written as
	\begin{equation*}
		\int_{0}^{R_0}e\left(\frac{u}{z}\right) \sum_{n=p}^{\infty}\frac{a_{n}}{\mu(n-1)}u^{n-1}\,du -\int_{R_0}^{\infty}e\left(\frac{u}{z}\right) \sum_{n=1}^{p-1}\frac{a_n}{\mu(n-1)}u^{n-1}\,du.
	\end{equation*}
	Then, we have
	\begin{equation}\label{eq.RemainderTwoSums}
	 \left|I_{\M,h}(\widehat{f}\,)(z)-\sum_{n=1}^{p-1}a_nz^n\right|\le I_{1,p}(z)+I_{2,p}(z),
	\end{equation}
	where
	\begin{align*}
		I_{1,p}(z)&=\left|\int_{0}^{R_0}e \left(\frac{u}{z}\right) \sum_{n=p}^{\infty}\frac{a_{n}}{\mu(n-1)}u^{n-1} \,du\right|,\\
		I_{2,p}(z)&=\left|\int_{R_0}^{\infty}e \left(\frac{u}{z}\right) \sum_{n=1}^{p-1}\frac{a_n}{\mu(n-1)}u^{n-1}\,du\right|.
	\end{align*}
	We first estimate $I_{1,p}(z)$.
	Since for every $u\in(0,R_0]$ we have $0<hu/h_1\le 1/2$, from~\eqref{eq.boundsCoeffBorelTransf} we get
	\begin{equation*}
		\sum_{n=p}^{\infty}\frac{|a_{n}|}{\mu(n-1)}u^{n-1}\le |\widehat{f}|_{\M,h} \frac{h}{h_1}\sum_{n=p}^{\infty}\left(\frac{hu}{h_1}\right)^{n-1}
	\le 2|\widehat{f}|_{\M,h}\left(\frac{h}{h_1}\right)^p u^{p-1}.
	\end{equation*}
	Hence,
\begin{equation}\label{eq.Estimates_f1}
I_{1,p}(z)\le 2|\widehat{f}|_{\M,h} \left(\frac{h}{h_1}\right)^p \int_{0}^{R_0}\left|e\left(\frac{u}{z}\right)\right| u^{p-1}\,du.
	\end{equation}
	Regarding $I_{2,p}(z)$, for $u\ge R_0$ and $1\le n\le p-1$ we have $(u/R_0)^{n-1}\le (u/R_0)^{p-1}$, so $u^{n-1}\le R_0^{n-1}u^{p-1}/R_0^{p-1}$. Again by~(\ref{eq.boundsCoeffBorelTransf}), and taking into account that $R_0=h_1/(2h)$, we may write
	\begin{align*}
		\sum_{n=1}^{p-1}\frac{|a_n|}{\mu(n-1)}u^{n-1}&\le |\widehat{f}|_{\M,h}\frac{hu^{p-1}}{h_1R_0^{p-1}} \sum_{n=1}^{p-1}\left(\frac{hR_0}{h_1}\right)^{n-1}\\
		&\le 2|\widehat{f}|_{\M,h}\frac{h}{h_1R_0^{p-1}}u^{p-1}
=|\widehat{f}|_{\M,h}\left(\frac{2h}{h_1}\right)^pu^{p-1}.
	\end{align*}
	Then,
	\begin{equation*}
		I_{2,p}(z)\le |\widehat{f}|_{\M,h} \left(\frac{2h}{h_1}\right)^p \int_{R_0}^{\infty}\left|e \left(\frac{u}{z}\right) \right|u^{p-1}\,du,
	\end{equation*}
and together with~\eqref{eq.RemainderTwoSums} and~\eqref{eq.Estimates_f1} we deduce
$$
\left|I_{\M,h}(\widehat{f}\,)(z)-\sum_{n=1}^{p-1}a_nz^n\right|\le |\widehat{f}|_{\M,h} \left(\frac{2h}{h_1}\right)^p \int_{0}^{\infty}\left|e \left(\frac{u}{z}\right) \right|u^{p-1}\,du.
$$
We estimate the last integral using first the second inequality in~\eqref{eq.Bounds_e_sector}, then the first one, and the fact that $e(x)>0$ for $x>0$, so that for every $z\in S_{\ga}$ and every $u>0$ we have
	$$
	|e(u/z)|\le K_3h_{\M}\left(K_4\frac{|z|}{u}\right)\le
	\frac{K_3}{K_1}\,e\left(\frac{K_2u}{K_4|z|}\right).
	$$
A change of variable and the right-hand estimates in~\eqref{eq.boundsMmomentslog} lead to
	\begin{align*}
		\int_{0}^{\infty}\left|e \left(\frac{u}{z}\right) \right|u^{p-1}\,du
		&\le \int_{0}^{\infty}\frac{K_3}{K_1}\,e\left(\frac{K_2u}{K_4|z|}\right) u^{p-1}\,du\\
		&=\frac{K_3}{K_1}\left(\frac{K_4|z|}{K_2}\right)^{p}\mu(p-1)
		\le \frac{K_3}{K_1}\left(\frac{K_4h_2}{K_2}\right)^{p}M_{p}|z|^{p}.
	\end{align*}
So, for every $p\in\N_0$ we have
	$$
	\left|I_{\M,h}(\widehat{f}\,)(z)-\sum_{n=1}^{p-1}a_nz^n\right|\le \frac{K_3|\widehat{f}|_{\M,h}}{K_1} \left(\frac{2K_4h_2h}{K_2h_1}\right)^pM_p|z|^p,\quad z\in S_{\ga},
	$$
and $I_{\M,h}(\widehat{f}\,)$ admits $\widehat{f}-a_0$ as its uniform $\{\M\}$-asymptotic expansion in $S_{\ga}$. We consider the map $T_{\M,h}$ defined in $\C[[z]]_{\M,h}$ as
	\begin{equation*}
T_{\M,h}(\widehat{f}\,)(z):=I_{\M,h}(\widehat{f}\,)(z)+a_0,\qquad z\in S_{\ga},
	\end{equation*}
which is obviously linear. Moreover, if we set $c:=2K_4h_2/(K_2h_1)>0$,
	\begin{equation*}
\left|T_{\M,h}(\widehat{f}\,)(z)-\sum_{n=0}^{p-1}a_nz^n\right|= \left|I_{\M,h}(\widehat{f}\,)(z)-\sum_{n=1}^{p-1}a_nz^n\right|\le \frac{K_3|\widehat{f}|_{\M,h}}{K_1} \left(ch\right)^pM_p|z|^p,\quad z\in S_{\ga},
	\end{equation*}
	which proves that $T_{\M,h}(\widehat{f}\,)\in\widetilde{\mathcal{A}}^u_{\M,ch}(S_{\ga})$ and, according to~\eqref{eq.NormUniformAsymptFixedType},
	$$
	\|T_{\M,h}(\widehat{f}\,)(z)\|_{\M,ch,\overset{\sim}{u}}\le \frac{K_3}{K_1}|\widehat{f}|_{\M,h},\ \ \widehat{f}\in\C[[z]]_{\M,h},
	$$
so that the continuity of $T_{\M,h}\colon \C[[z]]_{\M,h}\to \widetilde{\mathcal{A}}^u_{\M,ch}(S_{\ga})$ is obtained.
	
	$(ii)\Rightarrow (iii)$ Immediate for any weight sequence $\M$.
	
	$(iii)\Rightarrow (iv)$ It follows from~\eqref{equaContentionSurjectIntervals}, again valid for any weight sequence.
\end{proof}

\section{New surjectivity results in the Beurling Case}\label{sect.Improvements.Beurling}

In this last section we collect some results on surjectivity and existence of right inverses for the asymptotic Borel maps for Beurling ultraholomorphic classes. We split the results into two subsections according to the condition imposed on the sequence, $\dc$ or $\sm$.

\subsection{Continuous right inverses under derivation closedness}

A first result on the length of the interval $S_{(\hM)}$ was already provided by J. Schmets and M. Valdivia~\cite[Th. 4.4 and 4.6]{SchmetsValdivia00}, and it can be rephrased as follows.
Note that the target space for the Borel mapping in their results is the sequence space whose elements contain the derivatives of a given function at 0, subject to the natural growth restrictions; we are considering instead the target space consisting of the formal power series of asymptotic expansion, or the formal Taylor series, with a suitable control on their coefficients. Both spaces are naturally isomorphic, and both approaches are clearly equivalent.
Here, $\lfloor x\rfloor$ stands for the greatest integer which is less than or equal to the real number $x$.

\begin{theorem}\label{th.SchmVald.ext.oper.Beurling.implies.gamma}
Let $\hM$ be a regular sequence and $r>0$. If there is a global extension operator $U_{\M}:\CC[[z]]_{(\M)}\to\mathcal{A}_{(\hM)}(S_{r})$, then $\gamma(\M)>\lfloor r\rfloor$.
\end{theorem}

Later on, and regarding strongly regular sequences, the aforementioned result of V. Thilliez~\cite[Cor. 3.4.1]{Thilliez03} showed that $(0,\ga(\M))\subset S_{(\hM)}$, and A. Debrouwere~\cite[Cor. 1.3]{DebrouwereBorelRittBeurlingClasses} has recently proved that surjectivity comes with global extension operators. It is worth noting that~\cite[Th. 1.2]{DebrouwereBorelRittBeurlingClasses} gives a complete solution to the Borel-Ritt problem in non-uniform Beurling classes (which are not treated in this paper) defined by strongly regular sequences.

Going back to results without assuming the condition $\mg$, we first mention that the hypotheses in Theorem~\ref{th.SchmVald.ext.oper.Beurling.implies.gamma} can be improved by using some techniques included in~\cite{SchmetsValdivia00}, and Proposition 4.3 therein, in the same line of ideas that inspired the proof of a similar statement~\cite[Theorem 4.14(i)]{JimenezSanzSchindlInjectSurject} in the Roumieu case. Note that we will exchange the existence of the extension operator into just the surjectivity of the Borel map, and that $\mathcal{A}_{(\hM)}(S_{r})\subset \widetilde{\mathcal{A}}^u_{(\M)}(S_{r})$, what again weakens the forthcoming assumption.

\begin{theorem}\label{th.SurjectivityUniformAsymp.plus.dc}
Let $\hM$ be a regular sequence. If $r>0$ is such that $\widetilde{\mathcal{B}}:\widetilde{\mathcal{A}}^u_{(\M)}(S_{r})\to \CC[[z]]_{(\M)}$ is surjective, then $\ga(\M)>\lfloor r \rfloor $.
\end{theorem}

We include the proof for the reader's convenience. First, for $r\in\N$ we need to introduce the space $\mathcal{N}_{r,(\M)}([0,\infty))$ consisting of the functions $f\in\mathcal{C}^{\infty}([0,\infty))$ such that:
\begin{enumerate}[(a)]
 \item $f^{(pr+j)}(0)=0$ for every $p\in\N_0$ and $j\in\{1,\dots,r-1\}$ (this condition is empty when $r=1$),
 \item for every $h>0$ one has
$$
\sup_{p\in\N_0,\, x\in[0,\infty)}\frac{|f^{(pr)}(x)|}{h^pM_p}<\infty.
$$
\end{enumerate}

We recall the following crucial result.

\begin{prop}[\cite{SchmetsValdivia00}, Prop.\ 4.3]\label{prop.4.3.Nr.SchmetsValdivia}
Let $r\in\N$ and $\M$ be a weight sequence. If the map $\mathcal{B}_r:\mathcal{N}_{r,(\M)}([0,\infty))\longrightarrow \C[[z]]_{(\widecheck{\M})}$ sending $f$ to the formal power series $\sum^\infty_{p=0} (f^{(pr)}(0)/p!) z^p$ is surjective, then the sequence $\m$ satisfies the condition $(\gamma_r)$ or, in other words, $\ga(\M)>r$.
\end{prop}

\begin{proof}[Proof of Theorem~\ref{th.SurjectivityUniformAsymp.plus.dc}]
If $r\in(0,1)$, then it suffices to observe that $\widetilde{S}^u_{(\M)}$ is not empty in order to conclude, by Lemma~\ref{lemma.SurjectImpliesgammaMpositive}, that $\ga(\M)>0=\lfloor r\rfloor$ (note that $\dc$ has not been used in this case).

Suppose now that $r\ge 1$ and put $r_0=\lfloor r\rfloor\in\N$.
We will show that $\mathcal{B}_{r_0}:\mathcal{N}_{r_0,(\M)}([0,\infty))\longrightarrow \C[[z]]_{(\widecheck{\M})}$ is surjective, and then $\ga(\M)>r_0=\lfloor r\rfloor$ by Proposition~\ref{prop.4.3.Nr.SchmetsValdivia}, as desired.

Given $\widehat{g}=\sum^\infty_{p=0} a_p z^p\in\C[[z]]_{(\widecheck{\M})}$, we write $b_p:=a_p p!$ for all $p\in\N_0$, and for every $h>0$ there exists $C_1>0$ such that
\begin{equation}\label{equaBoundsbpProofDC}
|b_p|\le C_1h^pp!\widecheck{M}_p=C_1h^pM_p,\quad p\in\N_0.
\end{equation}
Consider the formal power series $\widehat{f}=\sum_{p=0}^{\infty}(-1)^{pr_0}b_pz^p\in\C[[z]]_{(\M)}$. By hypothesis, there exists $\psi\in\widetilde{\mathcal{A}}^u_{(\hM)}(S_{r})$ such that $\widetilde{\mathcal{B}}(\psi)=\widehat{f}$, and so there exists $C_2>0$ such that for every $p\in\N_0$ one has
\begin{equation}\label{equaAsympExpanPsiProofDC}
\Big|\psi(z)-\sum_{k=0}^{p-1}(-1)^{kr_0}b_kz^k \Big|\le C_2h^pM_{p}|z|^p,\qquad z\in S_r.
\end{equation}
The function $\varphi:S_{r/r_0}\to\C$ given by $\varphi(w)=\psi(w^{-r_0})-b_0$, is well defined and holomorphic in $S_{r/r_0}\supseteq S_1$. Moreover, according to~\eqref{equaAsympExpanPsiProofDC} for $p=1$, for every $w\in S_1$ one has
\begin{equation}\label{equaBoundsVarphiProofDC}
\left|\frac{\varphi(w)}{w}\right|=\frac{1}{|w|}|\psi(w^{-r_0})-b_0|\le \frac{C_2hM_{1}}{|w|^{r_0+1}}.
\end{equation}
So, the function $f:\R\to\C$ given by
\begin{equation*}
f(t)=\frac{1}{2\pi i}\int_{1-\infty\,i}^{1+\infty\,i}e^{tu}
\frac{\varphi(u)}{u}\,du
\end{equation*}
is well defined and continuous on $\R$. By the classical Hankel formula for the reciprocal Gamma function, for every natural number $p\ge 2$ and every $t\in\R$ we may write
\begin{equation}\label{equaRemainderProofDC}
f(t)-\sum_{k=1}^{p-1}(-1)^{kr_0}b_k\frac{t^{kr_0}}{(kr_0)!}=
\frac{1}{2\pi i}\int_{1-\infty\,i}^{1+\infty\,i}e^{tu}
\left(\frac{\varphi(u)}{u}- \sum_{k=1}^{p-1}\frac{(-1)^{kr_0}b_k}{u^{kr_0+1}}\right)\,du.
\end{equation}
Since, again by~\eqref{equaAsympExpanPsiProofDC}, we have
\begin{equation}\label{equaRemainderPhiProofDC}
\left|\frac{\varphi(u)}{u}- \sum_{k=1}^{p-1}(-1)^{kr_0}b_k\frac{1}{u^{kr_0+1}}\right|=
\frac{1}{|u|}\left|\psi(u^{-r_0})- \sum_{k=0}^{p-1}(-1)^{kr_0}b_k(u^{-r_0})^k\right|\le
\frac{C_2h^pM_p}{|u|^{pr_0+1}}
\end{equation}
for every $u\in S_1$, we can apply Leibniz's theorem for parametric integrals and deduce that the function
$$
f(t)-\sum_{k=1}^{p-1}(-1)^{kr_0}b_k\frac{t^{kr_0}}{(kr_0)!}
$$
belongs to $\mathcal{C}^{pr_0-1}(\R)$. Moreover, all of its derivatives of order $m\le pr_0-1$ at $t=0$ vanish, see the proof of~\cite[Theorem 4.14(i)]{JimenezSanzSchindlInjectSurject}.

As $p$ is arbitrary, we have that $f\in\mathcal{C}^{\infty}(\R)$ and, moreover,
$$
f^{(m)}(0)=\begin{cases}
  (-1)^{pr_0}b_p&\textrm{if $m=pr_0$ for some $p\ge 1$},\\
  0&\textrm{otherwise.}
\end{cases}
$$
Finally, we define the function
$$
F(t)=b_0+f(-t),\quad t\ge 0.
$$
Obviously, $F\in\mathcal{C}^{\infty}([0,\infty))$ and
$F^{(pr_0)}(0)=b_p$, $p\in\N_0$; $F^{(m)}(0)=0$ otherwise.
In order to conclude, we estimate the derivatives of $F$ of order $pr_0$ for arbitrary $p\in\N_0$. For $p=0$ and $t\ge 0$, we take into account~\eqref{equaBoundsbpProofDC} and \eqref{equaBoundsVarphiProofDC} in order to obtain that
\begin{equation}\label{equaBoundFtProofDC}
|F^{(0)}(t)|\le |b_0|+\frac{1}{2\pi}\int_{-\infty}^{\infty}e^{-t}
\frac{C_2hM_{1}}{|1+yi|^{r_0+1}}\,dy\le
C_1+\frac{C_2hM_1}{2\pi}
\int_{-\infty}^{\infty}\frac{1}{(1+y^2)^{(r_0+1)/2}}\,dy,
\end{equation}
and so $F$ is bounded.
For $p\ge 1$ we may write formula~\eqref{equaRemainderProofDC}  evaluated at $-t$ as
\begin{equation*}
f(-t)-\sum_{k=1}^{p}b_k\frac{t^{kr_0}}{(kr_0)!}=
\frac{1}{2\pi i}\int_{1-\infty\,i}^{1+\infty\,i}e^{-tz}
\left(\frac{\varphi(z)}{z}- \sum_{k=1}^{p}\frac{(-1)^{kr_0}b_k}{z^{kr_0+1}}\right)\,dz.
\end{equation*}
Then,
\begin{align*}
F^{(pr_0)}(t)&=b_p+\left(f(-t)- \sum_{k=1}^{p}b_k\frac{t^{kr_0}}{(kr_0)!}\right)^{(pr_0)}(t)\\
&=
b_p+\frac{1}{2\pi i}\int_{1-\infty\,i}^{1+\infty\,i}e^{-tz}
(-z)^{pr_0}\left(\frac{\varphi(z)}{z}- \sum_{k=1}^{p}\frac{(-1)^{kr_0}b_k}{z^{kr_0+1}}\right)\,dz,
\end{align*}
and we may apply~\eqref{equaBoundsbpProofDC}, and~\eqref{equaRemainderPhiProofDC} in order to obtain
\begin{equation}\label{equaBoundDerivFtProofDC}
|F^{(pr_0)}(t)|\le C_1h^pM_p+\frac{C_2h^{p+1}M_{p+1}}{2\pi} \int_{-\infty}^{\infty}\frac{1}{(1+y^2)^{(r_0+1)/2}}\,dy.
\end{equation}
From \eqref{equaBoundFtProofDC} and~\eqref{equaBoundDerivFtProofDC}, and since $\M$ satisfies $\dc$, we deduce that there exist $C_3,C_4>0$ and $H>1$ such that for every $p\in\N_0$ one has
$$
|F^{(pr_0)}(t)|\le C_1h^pM_p+C_3(Hh)^pM_p\le C_4(Hh)^pM_p,\quad t\ge 0.
$$
Since $h$ is arbitrary and $H$ does not depend on it, we see $F\in\mathcal{N}_{r_0,(\M)}([0,\infty))$ and $\mathcal{B}_{r_0}(F)=\widehat{g}$, and so $\mathcal{B}_{r_0}$ is surjective.
\end{proof}

In a recent work, A. Debrouwere~\cite{momentsdebrouwere} characterized the surjectivity of the asymptotic Borel map in the right half-plane for regular sequences in the terms of the gamma index associated with this sequence.

\begin{theorem}[\cite{momentsdebrouwere}, Theorem 7.4]\label{th.ExtensOperHalfplane}
	Suppose $\hM$ is a regular sequence. The following are equivalent:
	\begin{enumerate}[(i)]
		\item The Borel map $\widetilde{\mathcal{B}}\colon \mathcal{A}_{(\hM)}(S_1)\to\CC[[z]]_{(\M)}$ is surjective.
		\item There exists a global extension operator $U_{\M}:\CC[[z]]_{(\M)}\to\mathcal{A}_{(\hM)}(S_{1})$.
		\item $\gamma(\M)>1$.
	\end{enumerate}
\end{theorem}

Note that the implication $(ii)\Rightarrow(iii)$ corresponds to Theorem~\ref{th.SchmVald.ext.oper.Beurling.implies.gamma} for $r=1$, and that $(iii)\Rightarrow(i)$ was obtained by V. Thilliez, as already mentioned. Also, the implication $(i)\Rightarrow(iii)$ is slightly weaker than our previous result applied for $r=1$. However, the full equivalence is a powerful result, as $(ii)$ is deduced from any of the other two conditions.

As it occurs in the Roumieu case, see~\cite{JimenezSanzSchindlSurjectDC}, this information can be taken into the case of Beurling classes in a general sector by applying general Laplace
and Borel integral transforms of order $\al>0$, which basically arise from the classical transforms (inverse of each other) combined with ramifications of exponent $\al$.
We sketch the information needed, as the details can be found in Sections 5.5 and 5.6 of~\cite{balserutx}.

For $0<\al<2$, to the Laplace kernel function
$$
e_{\al}(z):=\frac{1}{\al}z^{1/\al}\exp(-z^{1/\al}),\qquad z\in S_{\al},
$$
there corresponds the moment function
$$
m_{\al}(\lambda):=\int_{0}^{\infty}t^{\lambda-1}e_{\al}(t)dt=\Gamma(1+\al \lambda), \qquad\text{Re}(\lambda)\ge 0,
$$
and the Borel kernel function
$$E_{\al}(z):=\sum_{p=0}^\infty \frac{z^p}{m_{\al}(p)}=
\sum_{p=0}^\infty \frac{z^p}{\Gamma(1+\al p)},\qquad z\in\C,$$
which is the classical Mittag-Leffler function of order $\al$.

Given a function $f$ holomorphic in a sector $S=S(d,\beta)$ (for some $\beta>0$) with suitable growth, the \textit{$\al$-Laplace transform of $f$ in a direction $\tau$ in $S$} is defined as
\begin{equation*}
	 (\mathcal{L}_{\al,\tau}f)(z):=\int_0^{\infty(\tau)}e_{\al}(u/z)f(u)\frac{du}{u},\quad
	|\arg(z)-\tau|<\al\pi/2,\ |z|\textrm{ small enough},
\end{equation*}
where integration is along the half-line parameterized by $t\in(0,\infty)\mapsto te^{i\tau}$. The family $\{\mathcal{L}_{\al,\tau}f\}_{\tau\textrm{\,in\,}S}$ defines, by analytic continuation, a function $\mathcal{L}_{\al}f$ named the \textit{$\al$-Laplace transform} of $f$, which is holomorphic in a sectorial region (see~\cite{balserutx} for details) bisected by $d$ of opening $\pi(\beta+\al)$.

Now, let $S=S(d,\be,r):=\{z\in S(d,\be)\colon |z|<r\}$ be a sector with $\be>\al$, and $f:S\to \C$ be holomorphic in $S$ and continuous at 0 (that is, the limit of $f$ at 0 exists when $z$ tends to 0 in every proper subsector of $S$). For $\tau\in\R$ such that $|\tau-d|<(\be-\al)\pi/2$ we consider a path
$\delta_{\al}(\tau)$ in $S$ consisting of a segment from the origin to a point $z_0$ with $\arg(z_0)=\tau+\al(\pi+\varepsilon)/2$ (for some
suitably small $\varepsilon\in(0,\pi)$), then the circular arc $|z|=|z_0|$ from $z_0$ to
the point $z_1$ on the ray $\arg(z)=\tau-\al(\pi+\varepsilon)/2$ (traversed clockwise), and
finally the segment from $z_1$ to the origin. The \textit{$\al$-Borel transform of $f$ in direction $\tau$} is defined as
$$
(\mathcal{B}_{\al,\tau}f)(u):=\frac{-1}{2\pi i}\int_{\delta_{\al}(\tau)}E_{\al}(u/z)f(z)\frac{dz}{z},\quad
u\in S(\tau,\varepsilon_0), \quad \varepsilon_0\textrm{ small enough}.
$$
The family
$\{\mathcal{B}_{\al,\tau}f\}_{\tau}$ defines the \textit{$\al$-Borel transform} of $f$, holomorphic in the sector $S(d,\be-\al)$ and denoted by $\mathcal{B}_{\al}f$.

In case $\al\ge 2$, $\mathcal{L}_{\al}f$ and $\mathcal{B}_{\al}f$ are defined by combining suitable ramification operators with the previous ones, see again~\cite{balserutx}.

The formal $\al$-Laplace and $\al$-Borel transforms, defined from $\C[[z]]$ into $\C[[z]]$, are respectively given by
$$\widehat{\mathcal{L}}_{\al}\big(\sum_{p=0}^{\infty}a_{p}z^{p}\big):= \sum_{p=0}^{\infty}\Gamma(1+\al p)a_{p}z^{p},\qquad \widehat{\mathcal{B}}_{\al}\big(\sum_{p=0}^{\infty}a_{p}z^{p}\big):= \sum_{p=0}^{\infty}\frac{a_{p}}{\Gamma(1+\al p)}z^{p}.$$
Note that the formal Laplace and Borel transforms, $\widehat{\mathcal{L}}_{\al}$ and $\widehat{\mathcal{B}}_{\al}$, are topological isomorphisms between the space $\C[[z]]_{(\M)}$ and $\C[[z]]_{(\M\cdot\bGa_{\al})}$, respectively $\C[[z]]_{(\M/\bGa_{\al})}$, for an arbitrary sequence $\M$.

The following result for the Roumieu case appeared in \cite[Th. 3.5]{JimenezSanzSchindlSurjectDC}. The fact that the constants $C$ and $c$, appearing in the next items, do not depend on the value of $h$, makes the result valid also for the Beurling case in a straightforward way, and its proof is therefore omitted. We use the notation $\M\cdot\bGa_{\al}$, respectively $\M/\bGa_{\al}$, for the sequences which are termwise product, resp. quotient, of $\M$ and $\bGa_{\al}=(\Gamma(1+\al p))_p$.

\begin{theorem}
\label{teorrelacdesartransfBL}
	Suppose $\M$ is an arbitrary sequence, and $\al,\gamma>0$. Let $f\in\widetilde{\mathcal{A}}^{u}_{[\M]}(S_\gamma)$
	and $f\sim^u_{[\M]}\widehat{f}$. Then, the following hold:
	\begin{enumerate}[(i)]
		\item For every $\be$ with $0<\be<\gamma$ one has $$\mathcal{L}_{\al}f\in \widetilde{\mathcal{A}}^{u}_{[\M\cdot\bGa_{\al}]}(S_{\be+\al}) \quad\textrm{and}\quad  \mathcal{L}_{\al}f\sim^u_{[\M\cdot\bGa_{\al}]} \widehat{\mathcal{L}}_{\al}\widehat{f}.
		$$
		Moreover, there exist $C,c>0$, depending only on $\al$, $\be$ and $\gamma$, such that for every $h>0$ and every $f\in\widetilde{\mathcal{A}}^{u}_{\M,h}(S_\gamma)$ one has
		 $\|\mathcal{L}_{\al}f\|_{\M\cdot\bGa_{\al},ch,\overset{\sim}{u}}\le C\|f\|_{\M,h,\overset{\sim}{u}}$, and so the maps $\mathcal{L}_{\al}\colon \widetilde{\mathcal{A}}^{u}_{\M,h}(S_\gamma)\to \widetilde{\mathcal{A}}^{u}_{\M\cdot\bGa_{\al},ch}(S_{\be+\al})$ and $\mathcal{L}_{\al}\colon \widetilde{\mathcal{A}}^{u}_{[\M]}(S_\gamma)\to \widetilde{\mathcal{A}}^{u}_{[\M\cdot\bGa_{\al}]}(S_{\be+\al})$ are continuous.
		\item Suppose $\gamma>\al$. For every $\be$ with $\al<\be<\gamma$ one has
		$$\mathcal{B}_{\al}f\in \widetilde{\mathcal{A}}^{u}_{[\M/\bGa_{\al}]}(S_{\be-\al})\quad \textrm{and}\quad \mathcal{B}_{\al}f\sim^u_{[\M/\bGa_{\al}]} \widehat{\mathcal{B}}_{\al}\widehat{f}.
		$$
		Moreover, there exist $C,c>0$, depending only on $\al$, $\be$ and $\gamma$, such that for every $h>0$ and every $f\in\widetilde{\mathcal{A}}^{u}_{\M,h}(S_\gamma)$ one has
		 $\|\mathcal{B}_{\al}f\|_{\M/\bGa_{\al},ch,\overset{\sim}{u}}\le C\|f\|_{\M,h,\overset{\sim}{u}}$, and so the maps $\mathcal{B}_{\al}\colon \widetilde{\mathcal{A}}^{u}_{\M,h}(S_\gamma)\to \widetilde{\mathcal{A}}^{u}_{\M/\bGa_{\al},ch}(S_{\be-\al})$ and $\mathcal{B}_{\al}\colon \widetilde{\mathcal{A}}^{u}_{[\M]}(S_\gamma)\to \widetilde{\mathcal{A}}^{u}_{[\M/\bGa_{\al}]}(S_{\be-\al})$ are continuous.
	\end{enumerate}
\end{theorem}

The use of Laplace and Borel transforms of arbitrary positive order allows us to generalize Theorem \ref{th.ExtensOperHalfplane} for arbitrary sectors. The idea of the proof for the Roumieu case~\cite[Th. 4.2]{JimenezSanzSchindlSurjectDC} applies to the Beurling case, we include it for the sake of completeness. We note that this result may also be deduced from the results in~\cite{DebrouwereBorelRittBeurlingClasses} about classes with non-uniform asymptotics, but we think it is interesting to provide an argument contained in our framework. This procedure makes Theorem~\ref{th.SurjectivityUniformAsymp.plus.dc} necessary.

\begin{theorem}\label{th.GlobalExtOperRegSeq}
	Suppose $\hM$ is a regular sequence, and let $r>0$. Each of the following statements implies the next one:
	\begin{enumerate}[(i)]
		\item $r<\gamma(\M)$.
		\item There exists a global extension operator $U_{\M,r}:\CC[[z]]_{(\M)}\to\widetilde{\mathcal{A}}^u_{(\M)}(S_{r})$.
		\item The Borel map $\widetilde{\mathcal{B}}:\widetilde{\mathcal{A}}^u_{(\M)}(S_{r})\to \CC[[z]]_{(\M)}$ is surjective.
		\item $r\le \gamma(\M)$.
	\end{enumerate}
\end{theorem}

\begin{proof}
	(i)$\Rightarrow$(ii) We consider two cases:
	\begin{enumerate}[({a}.1)]
		\item Suppose $r>1$, and take a real number $r'$ with $r<r'<\gamma(\M)$. The sequence $\boldsymbol{P}_1:=\hM/\bGa_{r'}$ satisfies $\dc$ and thanks to~\eqref{equa.gamma.cociente}, $\gamma(\boldsymbol{P}_1)=\gamma(\M)-r'+1>1$. By Lemma~\ref{lemma.gammaMgreaterthan1}, there exists a weight sequence $\boldsymbol{P}$ such that $\boldsymbol{p}\simeq\m/\overline{\boldsymbol{g}}^{r'-1}$, satisfies $\dc$ and $\gamma(\boldsymbol{P})=\gamma(\M)+1-r'>1$. Since the classes associated with $\boldsymbol{P}$ and $\M/\bGa_{r'-1}$ agree, Theorem~\ref{th.ExtensOperHalfplane} provides an extension operator
		$$
		 U\colon\CC[[z]]_{\left(\M/\bGa_{r'-1}\right)}\to \mathcal{A}_{\left(\hM/\bGa_{r'-1}\right)}(S_{1}).
		$$
		By Proposition~\ref{propcotaderidesaasin-RoumieuBeurling}.(i), we have that $\mathcal{A}_{\left(\hM/\bGa_{r'-1}\right)}(S_{1})\hookrightarrow \widetilde{\mathcal{A}}^u_{\left(\M/\bGa_{r'-1}\right)}(S_{1})$, and therefore this induces an extension operator
		$$
		 \widetilde{U}\colon\CC[[z]]_{\left(\M/\bGa_{r'-1}\right)}\to \widetilde{\mathcal{A}}^u_{\left(\M/\bGa_{r'-1}\right)}(S_{1}).
		$$
		Since $r<r'=1+(r'-1)$, Theorem~\ref{teorrelacdesartransfBL}.(i) implies that the composition
		$$
		U_{\M,r}:=\mathcal{L}_{r'-1}\circ\widetilde{U}\circ \widehat{\mathcal{B}}_{r'-1}:\CC[[z]]_{(\M)}\to \widetilde{\mathcal{A}}^u_{(\M)}(S_{r})		$$
		is the extension operator we were looking for.
		\item If $r\le 1$, consider $\al$ such that $\al+r>1$, and take $r'$ with $r<r'<\gamma(\M)$. The weight sequence $\M\cdot{\bGa}_{\al}$ satisfies $\dc$ and $\gamma(\M\cdot{\bGa}_{\al})=\gamma(\M)+\al>r'+\al>1$. By item (a.1), there exists an extension operator
		$$
		U\colon\CC[[z]]_{(\M\cdot{\bGa}_{\al})}\to \widetilde{\mathcal{A}}^u_{(\M\cdot{\bGa}_{\al})}(S_{r'+\al}).
		$$
		Now, Theorem~\ref{teorrelacdesartransfBL}.(ii) implies that
		$$
		U_{\M,r}:=\mathcal{B}_{\al}\circ U\circ \widehat{\mathcal{L}}_{\al}:\CC[[z]]_{(\M)}\to \widetilde{\mathcal{A}}^u_{(\M)}(S_{r})
		$$
is the desired extension operator.
	\end{enumerate}\par\noindent
	(ii)$\Rightarrow$(iii)  The existence of $U_{\M,r}$ implies that the corresponding Borel map $ \widetilde{\mathcal{B}}:\widetilde{\mathcal{A}}^u_{(\M)}(S_{r})\to\CC[[z]]_{(\M)}$ is surjective in $S_r$.\par\noindent
	(iii)$\Rightarrow$(iv) Let us see that $r\leq \gamma(\M)$. We again have different cases:
\begin{enumerate}[({b}.1)]
	\item If $0<r<1$, consider positive real numbers $\al,r'$ with $1-\al<r'<r$. By applying the Laplace transform $\mathcal{L}_{\al}\colon \widetilde{\mathcal{A}}^{u}_{(\M)}(S_{r})\to \widetilde{\mathcal{A}}^{u}_{(\M\cdot\bGa_{\al})}(S_{r'+\al})$, Theorem~\ref{teorrelacdesartransfBL}.(i) shows that the map
	$$
	 \widetilde{\mathcal{A}}^u_{(\M\cdot\bGa_{\al})}(S_{r'+\al})\to\CC[[z]]_{(\M\cdot\bGa_{\al})},
	$$
	 is surjective. Observe that $r'+\al>1$, so we deduce by restriction to the half-plane $S_1$ that, according to Proposition~\ref{propcotaderidesaasin-RoumieuBeurling}.(ii), also the map
	 $$
	 	 \mathcal{A}_{(\hM\cdot\bGa_{\al})}(S_1)\to\CC[[z]]_{(\M\cdot\bGa_{\al})},
	 $$
	 is surjective. Theorem~\ref{th.ExtensOperHalfplane} implies then that $\gamma(\M\cdot\bGa_{\al})>1$ or, equivalently by~\eqref{equa.gamma.producto}, $\gamma(\M)>1-\al$. Since $\al$ can be chosen arbitrarily while keeping $1-\al<r$, we deduce $\gamma(\M)\ge r$.
	\item If $r\in\N$, we know that $\gamma(\M)>r$ by Theorem~\ref{th.SurjectivityUniformAsymp.plus.dc}.
	\item If $r\in(1,\infty)\setminus\N$, again by Theorem~\ref{th.SurjectivityUniformAsymp.plus.dc} we deduce that $\gamma(\M)>\lfloor r\rfloor$, so that the sequence $\boldsymbol{P}_1:=\hM/\bGa_{\lfloor r\rfloor}$ is such that $\gamma(\boldsymbol{P}_1)>1$ by using the properties of gamma index. 
	Hence, by Lemma~\ref{lemma.gammaMgreaterthan1} there exists a weight sequence $\boldsymbol{P}$ such that $\boldsymbol{P}\approx\M/\bGa_{\lfloor r\rfloor}$, $\gamma(\boldsymbol{P})=\gamma(\M)-{\lfloor r\rfloor}$ and  $\boldsymbol{P}$ will also satisfy $\dc$. Consider a value $r'$ with $\lfloor r\rfloor<r'<r$. By applying the Borel transform $\mathcal{B}_{\lfloor r\rfloor}\colon \widetilde{\mathcal{A}}^{u}_{(\M)}(S_{r})\to \widetilde{\mathcal{A}}^{u}_{(\M/\bGa_{\lfloor r\rfloor})} (S_{r'-\lfloor r\rfloor})$, Theorem~\ref{teorrelacdesartransfBL}.(ii) shows that the map
	$$
	\widetilde{\mathcal{A}}^{u}_{(\M/\bGa_{\lfloor r\rfloor})} (S_{r'-\lfloor r\rfloor})\to\CC[[z]]_{(\M/\bGa_{\lfloor r\rfloor})},
	$$
	is surjective, or equivalently, thanks to the equivalence $\boldsymbol{P}\approx\M/\bGa_{\lfloor r\rfloor}$, the map
	$$
	\widetilde{\mathcal{A}}^{u}_{(\boldsymbol{P})} (S_{r'-\lfloor r\rfloor})\to\CC[[z]]_{(\boldsymbol{P})},
	$$
	is also surjective. Since $ r'-\lfloor r\rfloor\in(0,1)$, we may invoke item (b.1) and deduce that $\gamma(\boldsymbol{P})\ge r'-\lfloor r\rfloor$, what amounts to $\gamma(\M)\ge r'$. We conclude by making $r'$ tend to $r$.
\end{enumerate}
\end{proof}

\subsection{Surjectivity for Beurling classes under condition $\sm$}

We end by proving a surjectivity result for Beurling classes when their defining weight sequence satisfies the new condition $\sm$. The technique used by V. Thilliez in~\cite[Th. 3.4.1]{Thilliez03} will be followed, and we first need to recall two auxiliary results from the work of J. Chaumat and A.-M. Chollet~\cite{chaucho}.

\begin{lemma}[\cite{chaucho}, Lemme 14]\label{lemma-ChauCho-Beurlingsucesiones} Let $\LL=(L_p)_p$ be a sequence of nonnegative real numbers and $\M=(M_p)_p$ be a sequence of positive real numbers. The following conditions are equivalent:
	\begin{enumerate}[(i)]
		\item For all $h>0$, there exists a constant $C(h)>0$ such that $L_p\leq C(h) h^p M_p$ for every $p\in\N_0$.
		\item $\displaystyle\lim_{p\to\infty}\left(\frac{L_p}{M_p}\right)^{1/p}=0$.
		\item There exists a sequence $\bbepsilon=(\varepsilon_p)_{p\in\N_0}$ of positive real numbers tending to zero such that $L_p\leq \varepsilon_0\varepsilon_1\cdots \varepsilon_{p-1} M_p$, $p\in\N$.
	\end{enumerate}
	Moreover, if (i) or (ii) are satisfied, then (iii) holds true for a nonincreasing sequence $\bbepsilon$.
\end{lemma}

\begin{lemma}[\cite{chaucho}, Lemme 16]\label{lemma-Chaucho-3sucesiones}
Let $\boldsymbol{A}=(A_p)_p$ be a sequence of nonnegative real numbers such that $\sum_{p=0}^{\infty} A_p$ is convergent, and let $\boldsymbol{B}=(B_p)_p$ and $\boldsymbol{D}=(D_p)_p$ be sequences of positive real numbers such that $\lim_{p\to\infty}B_p=0$, and $\boldsymbol{D}$ is nonincreasing and $\lim_{p\to\infty}D_p=0$. Then, there exists a nondecreasing sequence $\boldsymbol{E}=(E_p)_p$ of positive real numbers such that:
	\begin{enumerate}[(i)]
		\item $\lim_{p\to\infty}E_p=\infty$.
		\item $
\sum_{p=q}^{\infty}E_pA_p\leq 8 E_q \sum_{p=q}^{\infty}A_p$, $q\in\N_0$.
		\item The sequence $\boldsymbol{E}\cdot\boldsymbol{D}=(E_pD_p)_p$ is nonincreasing.
		\item $\lim_{p\to\infty}E_pB_p=0.$
	\end{enumerate}
\end{lemma}

The next result is an adaptation of a similar result,~\cite[Prop. 17]{chaucho}, in which the condition $\mg$ has now been substituted by $\sm$.

\begin{theorem}\label{theorem-Chaucho-BeurlingRoumieu} Let $\LL=(L_p)_p$ be a weight sequence satisfying $\snq$ and $\sm$.
If $\boldsymbol{A}=(A_p)_p$ is a sequence of nonnegative real numbers such that for all $h>0$ there exists $C(h)>0$ such that
$A_p\leq C(h) h^p L_p$ for every $p\in\N_0$, then there exists a weight sequence $\boldsymbol{K}=(K_p)_p$ which satisfies $\snq$, $\sm$ and such that:
\begin{enumerate}[i)]
	\item There exists a constant $D>0$ such that $A_p \leq D K_p$, for all $p\in\N_0$.
	\item For all $h>0$, there exists $C'(h)>0$ such that
	$K_p\leq C'(h) h^p L_p$, $p\in\N_0$. 	
\end{enumerate}
\end{theorem}

\begin{proof} By Lemma \ref{lemma-ChauCho-Beurlingsucesiones}, there exists a nonincreasing sequence $\bbepsilon=(\varepsilon_p)_{p\in\N_0}$ which tends to zero, and such that
\begin{equation}\label{Desigualdad-LpMp}
	A_p\leq \varepsilon_0\varepsilon_1\cdots \varepsilon_{p-1} L_p,\qquad p\in\N.
\end{equation}
Consider the sequence $(u_p)_{p\in\N_0}$ defined as
\begin{equation*}
	u_p=\frac{1}{(p+1)\ell_p},
\end{equation*}
where $(\ell_p)_{\in\N_0}$ is the sequence of quotients of $\LL$.
$(u_p)_{p\in\N_0}$ is nonincreasing and tends to zero, and since $\LL$ satisfies $\snq$, there exists some constant $A>0$ with
\begin{equation}\label{Chaucho-M-snq}
	\sum_{p=q}^{\infty} u_p\leq A (q+1)u_q, \qquad q\in\N_0.
\end{equation}
As $((p+1)u_p)_{p\in\N_0}$ is nonincreasing and tends zero, we can apply Lemma \ref{lemma-Chaucho-3sucesiones} with
$$
A_p=u_p,\qquad B_p=\max\{\varepsilon_p,(p+1)u_p\},\qquad D_p=(p+1)u_p,\qquad p\in\N_0.
$$
So, there exists a nondecreasing sequence $\boldsymbol{E}$ which tends to $\infty$, and such that:
\begin{align}
& \sum_{p=q}^{\infty}u_pE_p\leq 8 E_q \sum_{p=q}^{\infty}u_p, \qquad q\in\N_0.\label{Chaucho-consequence1}\\
& \text{The sequence $((p+1)u_pE_p)_{p\in\N_0}$ is nonincreasing.}\label{Chaucho-consequence2}\\
& \lim_{p\to\infty}\varepsilon_pE_p=0,\qquad \lim_{p\to\infty}(p+1)u_pE_p=0.\label{Chaucho-consequence3}
\end{align}
Let us consider the sequence $(k_p)_{p\in\N_0}$ defined as
$$
k_p=\frac{\ell_p}{E_p}=\frac{1}{(p+1)u_pE_p}, \qquad p\in\N_0.
$$
Then, from \eqref{Chaucho-M-snq} and \eqref{Chaucho-consequence1} we deduce that
$$
\sum_{p=q}^{\infty} \frac{1}{(p+1)k_p}=\sum_{p=q}^{\infty} u_pE_p\leq 8E_q \sum_{p=q}^{\infty} u_p\leq 8A(q+1)u_qE_q=8A\frac{1}{k_q},\qquad q\in\N_0.
$$
Therefore, the sequence $\boldsymbol{K}$ defined as
$$
K_0=1,\qquad K_p=k_0k_1\cdots k_{p-1}, \qquad p\in\N,
$$
satisfies $\snq$. Moreover, the sequence $\boldsymbol{K}$ is a weight sequence due to \eqref{Chaucho-consequence2} and \eqref{Chaucho-consequence3}.
Now, since $\boldsymbol{E}$ is nondecreasing we deduce that
$$
\frac{k_{p+1}}{k_p}=\frac{\ell_{p+1}E_p}{\ell_pE_{p+1}}\leq \frac{\ell_{p+1}}{\ell_p},\qquad p\in\N_0,
$$
and so $\boldsymbol{K}$ satisfies $\sm$ too (with the same constants as $\LL$).
Taking into account that
\begin{equation}\label{Chaucho-igualdadKpMp}
	K_p=\frac{1}{E_0\dots E_{p-1}}L_p, \qquad p\in\N,
\end{equation}
we deduce from \eqref{Desigualdad-LpMp} that
\begin{equation*}
	A_p\leq \varepsilon_0\varepsilon_1\cdots \varepsilon_{p-1} L_p=\varepsilon_0E_0\varepsilon_1E_1\cdots \varepsilon_{p-1}E_{p-1}K_p,\qquad p\in\N.
\end{equation*}
We observe in \eqref{Chaucho-consequence3} that the sequence $(\varepsilon_pE_p)_{p\in\N_0}$ tends to zero and, therefore, Lemma \ref{lemma-ChauCho-Beurlingsucesiones} provides for all $t>0$ a constant $D(t)>0$ such that
$A_p\leq D(t)t^pK_p$ for every $p\in \N_0$ and, in particular, for $t=1$ we obtain that $A_p\leq D(1)K_p$.
Finally, from \eqref{Chaucho-igualdadKpMp} and the fact that $(1/E_p)_{p\in\N_0}$ tends to zero, we can apply again Lemma \ref{lemma-ChauCho-Beurlingsucesiones} and we deduce that for all $h>0$ there exists some constant $C'(h)>0$ such that
$K_p\leq C'(h)h^pL_p$ for every $p\in \N_0$.
\end{proof}

We are ready for the proof of our last result.

\begin{theorem}\label{Thilliez-Beurlingcase}
Let $\M$ be a weight sequence with $\ga(\M)>0$ and that satisfies $\sm$, and let $0<r<\ga(\M)$ be given.
	Then, the Borel map $\widetilde{\mathcal{B}}\colon \widetilde{\mathcal{A}}^u_{(\M)}(S_{r})\to\C[[z]]_{(\M)}$ is surjective, and so $(0,\ga(\M))\subset\widetilde{S}_{(\M)}^u$.
\end{theorem}

\begin{proof}
Let $\widehat{f}=\sum_{p=0}^\infty a_pz^p\in\C[[z]]_{(\M)}$ be given. By the definition of $\C[[z]]_{(\M)}$ (see also~\eqref{eq.defBanachFormalPowerSeries}), for every $h>0$ there exists $C(h)>0$ such that
	\begin{equation}\label{eq.BeurlingBoundsCoeff}
		|a_p|\le C(h)h^{p}M_{p},\quad p\in\N_0.
	\end{equation}
On the one hand, by the properties of the gamma index we have that $\gamma(\M^{1/r})=\ga(\M)/r>1$, and Lemma~\ref{lemma.gammaMgreaterthan1} provides the existence of a weight sequence $\LL=(L_p)_p$ such that
$\gamma(\widehat{\LL})>1$ and of a constant $a>0$ such that $a^{-1}(\ell_p(p+1))\leq m_p^{1/r}\leq a(\ell_p(p+1))$ for all $p\in\N_0$.
In particular, we have that $\gamma(\LL)>0$, and so $\LL$ satisfies $\snq$. Moreover, it is clear that $a^{-p}p!L_p\le M_p^{1/r}\le a^pp!L_p$ for every $p$, and so $\M\approx (\widehat{\LL})^r$, what implies that the classes defined by these sequences coincide. Because of the stability properties of $\sm$ described in Subsection~\ref{subsectstrregseq}, $\LL$ inherits $\sm$ from $\M$.

On the other hand, from~\eqref{eq.BeurlingBoundsCoeff} we obtain
$$
\frac{|a_p|^{1/r}}{p!}\le C(h)^{1/r}(ah^{1/r})^{p}L_{p},\ \ p\in\N_0,
$$
and so we are in a position to apply Theorem~\ref{theorem-Chaucho-BeurlingRoumieu} to the sequences $\boldsymbol{A}=(|a_p|^{1/r}/p!)_p$ and $\LL$. Hence, there exists a weight sequence $\boldsymbol{K}=(K_p)_p$ which satisfies $\snq$ (i. e., $\ga(\boldsymbol{K})>0$) and $\sm$, such that there exists $D>0$ with $|a_p|^{1/r}/p! \leq D K_p$ for all $p\in\N_0$, and such that for all $h>0$, there exists $C'(h)>0$ such that
\begin{equation}\label{eq.K.strongly.less.L}
K_p\leq C'(h) h^p L_p,\ \  p\in\N_0.
\end{equation}
The first estimates state that $|a_p|\le D^r(p!K_p)^r$, and so $\widehat{f}\in\C[[z]]_{\{\NN\}}$ for the weight sequence $\NN:=(\widehat{\boldsymbol{K}})^r$. Again $\NN$ inherits $\sm$ from $\boldsymbol{K}$, and moreover
$\ga(\NN)=r(\ga(\boldsymbol{K})+1)>r$. So,
we can apply Theorem~\ref{Theorem-surjectivity-logcondition} to deduce that $\widetilde{\mathcal{B}}\colon \widetilde{\mathcal{A}}^u_{\{\NN\}}(S_{r})\to\C[[z]]_{\{\NN\}}$ is surjective. Hence, there exists $f\in\widetilde{\mathcal{A}}^u_{\{\NN\}}(S_{r})$ such that $\widetilde{\mathcal{B}}(f)=\widehat{f}$. Finally, observe that from~\eqref{eq.K.strongly.less.L} we get
$$
N_p=(p!K_p)^r\le C'(h)^r(h^r)^p(\widehat{L}_p)^r,\ \ p\in\N_0,
$$
and so $\widetilde{\mathcal{A}}^u_{\{\NN\}}(S_{r})\subset
\widetilde{\mathcal{A}}^u_{((\widehat{\LL})^r)}(S_{r})= \widetilde{\mathcal{A}}^u_{(\M)}(S_{r})$, from where the conclusion follows.
\end{proof}

\vskip.2cm
\noindent \textbf{Funding} \ The first three authors are partially supported by the Spanish Ministry of Science and Innovation under the project PID2022-139631NB-I00. The research of the fourth author was funded in whole or in part by the Austrian Science Fund (FWF)  10.55776/PAT9445424.\par

%
%
%
\vskip.2cm
\noindent\textbf{Affiliations}:\\
\noindent Javier~Jim\'{e}nez-Garrido:\\
Departamento de Matem\'aticas, Estad{\'\i}stica y Computaci\'on\\
Universidad de Cantabria\\
Avda. de los Castros, s/n, 39005 Santander, Spain\\
Instituto de Investigaci\'on en Matem\'aticas IMUVA, Universidad de Va\-lla\-do\-lid\\
ORCID: 0000-0003-3579-486X\\
E-mail: jesusjavier.jimenez@unican.es\\

\vskip.1cm
\noindent
Ignacio Miguel-Cantero:\\
Instituto de Investigaci\'on en Matem\'aticas IMUVA\\
Universidad de Va\-lla\-do\-lid\\
ORCID: 0000-0001-5270-0971\\
E-mail: ignacio.miguel@uva.es\\

\vskip.1cm
\noindent Javier~Sanz:\\
Departamento de \'Algebra, An\'alisis Matem\'atico, Geometr{\'\i}a y Topolog{\'\i}a\\
Universidad de Va\-lla\-do\-lid\\
Facultad de Ciencias, Paseo de Bel\'en 7, 47011 Valladolid, Spain.\\
Instituto de Investigaci\'on en Matem\'aticas IMUVA\\
ORCID: 0000-0001-7338-4971\\
E-mail: javier.sanz.gil@uva.es\\

\vskip.1cm
\noindent Gerhard~Schindl:\\
Fakult\"at f\"ur Mathematik, Universit\"at Wien,
Oskar-Morgenstern-Platz~1, A-1090 Wien, Austria.\\
ORCID: 0000-0003-2192-9110\\
E-mail: gerhard.schindl@univie.ac.at
\end{document}